\documentclass[12pt]{amsart}

\def\urlfont{\DeclareFontFamily{OT1}{cmtt}{\hyphenchar\font='057}
              \normalfont\ttfamily \hyphenpenalty=10000}
\def\pplogo{\vbox{\kern-\headheight\kern -15pt
\halign{##&##\hfil\cr&{
\ppnumber}\cr\rule{0pt}{2.5ex}&\ppdate\cr}
}}
\makeatletter
\def\ps@firstpage{\ps@empty \def\@oddhead{\hss\pplogo}%
  \let\@evenhead\@oddhead 
}
\def\maketitle{\par
 \begingroup
 \def\thefootnote{\fnsymbol{footnote}}
 \def\@makefnmark{\hbox
 to 0pt{$^{\@thefnmark}$\hss}}
 \if@twocolumn
 \twocolumn[\@maketitle]
 \else \newpage
 \global\@topnum\z@ \@maketitle \fi\thispagestyle{firstpage}\@thanks
 \endgroup
 \setcounter{footnote}{0}
 \let\maketitle\relax
 \let\@maketitle\relax
 \gdef\@thanks{}\gdef\@author{}\gdef\@title{}\let\thanks\relax}
\renewcommand{\thetable}{\@Alph\c@table}
\makeatother

\def\ppnumber{\vbox{\baselineskip14pt
\hbox{DUKE-CGTP-00-05}
\hbox{IASSNS-HEP-00/27}
\hbox{math.AG/0005196}
}}
\def\ppdate{Revised May, 2000}
\date{} 

\title[Euler characteristic of threefolds]{Group representations and the
Euler characteristic of elliptically fibered \\Calabi--Yau threefolds}

\author[A. Grassi and D. R. Morrison]{Antonella Grassi and David
R. Morrison}

\address{Department of Mathematics, University of Pennsylvania,
Philadelphia, PA 19104}
\email{grassi@math.upenn.edu}

\address{School of Natural Sciences, Institute for Advanced Study,
Princeton, NJ 08540, and\\
Department of Mathematics, Duke University, Durham, NC
27708-0320}
\email{drm@math.duke.edu}

\thanks{Research partially supported by the Harmon Duncombe foundation,
by the Institute for Advanced Study, and
by National Science Foundation grants DMS-9401447, DMS-9401495,
DMS-9627351 and DMS-9706707. We thank the Institute for Advanced Study,
 the Mathematisches Forschunginstitut  Oberwolfach, and
the Institute for Theoretical Physics, Santa Barbara,
for hospitality during  various stages of this project.}

\theoremstyle{plain}
\newtheorem{theorem}{Theorem}
\newtheorem{corollary}[theorem]{Corollary}
\newtheorem{proposition}[theorem]{Proposition}
\newtheorem{lemma}[theorem]{Lemma}
\newtheorem{definition}[theorem]{Definition}
\numberwithin{theorem}{section}
\numberwithin{equation}{section}

\theoremstyle{remark}
\newtheorem{remark}[theorem]{Remark}
\newcommand{\Tr}{\operatorname{Tr}}
\newcommand{\tr}{\operatorname{tr}}
\newcommand{\adj}{\operatorname{adj}}
\newcommand{\fund}{\operatorname{fund}}
\newcommand{\fundbar}{\mathop{\overline{\fund}}}
\newcommand{\spinrep}{\operatorname{spin}}
\newcommand{\vect}{\operatorname{vect}}
\newcommand{\ad}{\operatorname{ad}}
\newcommand{\End}{\operatorname{End}}
\newcommand{\trace}{\operatorname{trace}}

\def\bsa{\boldsymbol\Sigma}
\def\bs1{\boldsymbol\Sigma_1}
\def\bsi{\boldsymbol\Sigma_i}
\def\bsj{\boldsymbol\Sigma_j}
\def\bso{\boldsymbol\Sigma_0}

\def \a{\alpha }

\def \g{\gamma}
\def\so{{S _0}}
\def\pl{\mathbb P^1}
\def\bpt{\mathbb P^1 \times \mathbb P^1}
\def\p2{\mathbb P^2}
\def\r{\boldmath  $\mathcal R$}
\def\rloc#1{$\mathcal R_{#1}$}

\def\c{\chi _{top}}
\def\kb{ K_B}
\def\ks{K_B \cdot \bs1}
\def\cc{{\bs1}^2}
\def\cci{{\bsi}^2}
\def\su{\operatorname{SU}}
\def\spin{\operatorname{Spin}}
\def\sp{\operatorname{Sp}}
\def\dim{\operatorname{dim}}
\def\rk{\operatorname{rk}}
\def\so{\operatorname{SO}}
\def\GL{\operatorname{GL}}
\def\f{\frac{1}{2}}
\def\Z{\mathbb Z}
\newcommand{\nonmin}{NM}
\newcommand{\terminal}{NSR}

\renewcommand{\div}{\operatorname{div}}

\begin{document}

\begin{abstract}
To every elliptic Calabi--Yau threefold with a section $X$ there can
be associated 
a Lie group $G$ and a representation $\rho$ of that group.  The group
is determined from the Weierstrass model, which has singularities
that are generically rational double points; these double points
lead to local factors of $G$ which are either the corresponding
A-D-E groups or some associated non-simply laced groups.
  The representation $\rho$ is a 
sum of representations coming from the local factors of $G$,
and of other representations which can be associated to the
points at which the singularities are worse than generic.

This construction first arose in physics, and the requirement of
anomaly cancellation in the associated physical theory makes some surprising
predictions about the connection between $X$ and $\rho$.
In particular, an explicit formula (in terms of $\rho$) for the Euler characteristic
of $X$ is predicted.  We give a purely mathematical proof 
of that formula in this paper, introducing along the way a new invariant
of elliptic Calabi--Yau threefolds.  We also verify
 the other geometric predictions which are consequences
of anomaly cancellation, under some (mild) hypotheses about
the types of singularities which occur.

 As a  byproduct we also discover a novel relation between
the Coxeter number and the rank in the case of the simply laced groups in the
``exceptional series'' studied by Deligne.
\end{abstract}

\maketitle

 It was noted by Du Val \cite{DV} that certain surface singularities,
now known as rational double points,
are classified by the Dynkin diagrams of the
simply laced Lie groups\footnote{More precisely, Du Val recognized the
combinatorial structure as occurring in the theory of finite reflection
groups; the connection to Lie groups was made soon thereafter by Coxeter
\cite{coxeter:annals,coxeter:weyl}.} of type $A_n, D_n, E_6, E_7, E_8$.
Du Val pointed out that the Dynkin diagram is the dual diagram to the
intersection configuration
of the exceptional divisors in the minimal resolution of the singularities.
 Further connections between these singularities and Lie groups were
subsequently discovered by Brieskorn and Grothendieck
\cite{Bri:Nice}.

 The resolutions of rational double points are crepant, that is,
the pullback of the canonical divisor on the singular variety is
the canonical divisor on the smooth minimal resolution.
 In particular, if the singular variety has trivial canonical class,
so does its desingularization.

One characterization of rational double points is as quotients of ${\mathbb
C}^2$ by  finite
subgroups of $SL(2,{\mathbb C})$ \cite{DV:book}.
  Much recent work has been done by looking at
the quotient of ${\mathbb C}^3$ by a finite
subgroup of $SL(3,{\mathbb C})$  (see for example \cite{IR,IN,R}).
In this paper we consider another natural generalization of
the above set up.

It turns out that the singularities of the Weierstrass model of an
elliptic surface are also rational double points. If the singular surface
satisfies the Calabi--Yau condition, so does its resolution; the
Calabi--Yau condition can be expressed as a condition on the base
of the elliptic fibration and the discriminant locus.
 Furthermore the ranks of the A-D-E groups contribute to the rank of the
Picard group of the minimal resolution of the Weierstrass model,
as well as its topological Euler characteristic.

In this paper we investigate a similar situation one dimension
higher, namely, elliptic Calabi--Yau threefolds which
are resolutions of Weierstrass models.
 Here the singularities are only generically rational double points,
yet it is possible to associate a group $G$ to the singularities,
obtaining all the Dynkin diagrams (including the non-simply laced ones).
  We restrict our attention to Weierstrass models with a
minimal
resolution which is a flat elliptic fibration satisfying the Calabi--Yau
condition---the existence of a flat resolution excludes some non-generic
singularities.

The threefold also determines a specific representation of $G$
(known in the physics literature as the ``matter representation'')
whose irreducible summands can be
described in terms of the degenerations
of the general singularity; conversely,
 once one chooses the representations which might occur, the geometry of
the Calabi--Yau
is completely determined by  some relations in representation theory.

 As a  byproduct of our analysis,
we also discover a novel relation between
the Coxeter number and the rank in the case of the simply laced groups in the
``exceptional series'' studied by Deligne \cite{deligne}.

\medskip

 This work was first motivated by
the problem of verifying in this context the vanishing
of an ``anomaly''
coming from string theory (see Section \ref{sec:physics})
 and by completing  a dictionary between the geometry of the
Calabi--Yau and the corresponding quantities in quantum field theory. We can
in fact interpret the vanishing of the anomaly  as a formula for the
Euler characteristic
 of the Calabi--Yau manifold,  a formula which
was quite unexpected.

 We first formally define an invariant {\r} (see Section
\ref{sec:firstlook})
and show how certain representations of $G$
appear in {\r} when the ``general'' double point degenerates
to a worse singularity. We will then show how the geometry of the Calabi--Yau
and its degenerations are naturally, yet surprisingly,  related to the same
representations
occurring in {\r}  (see Section \ref{sec:sadov}).
{}From the string theory point of view this is explained by considering
a quantum field theory associated to $X$, which suffers from potential gauge
and gravitational anomalies.  Some of these anomalies can be ``cancelled'' by
an analogue of the Green--Schwarz mechanism, while others (which occur
as certain coefficients in a formal expression in the curvature) are
required to vanish identically.  The vanishing of the latter leads to the
formula for {\r}, while the existence of the Green--Schwarz mechanism imposes
the other geometric constraints.  Note that
our arguments and definitions, while inspired by  the physics of string
theory, are in the realm of mathematics only; the explicit dictionary between
mathematical and field-theoretic
quantities is developed in \cite{AD}.

\medskip

In Section \ref{sec:G} we discuss how we can associate a group $G$ to an
elliptic
threefold which is a resolution of a Weierstrass model; in Section
  \ref{sec:firstlook}  we introduce
the invariants {\r} and $H_{ch}$ and we present some first properties
of {\r}.
  In the following Section \ref{sec:physics} we sketch some background from
physics.
 Section \ref{sec:charged} shows that the results
of our Main Theorem \ref{th:main} agree with some  predictions
from the physics literature.

 After stating the working assumptions and some notation  in Section
\ref{sec:wa},
we describe an algorithm to compute the fundamental invariant {\boldmath
$\mathcal R$}
from the singularities of the Weierstrass model (Section \ref{sec:ec}).
  We show that the group $G$ determines most of the terms occurring in
{\boldmath  $\mathcal R$}
and we present these in  Appendix I.

 The other terms come from the degeneration of the
``generic'' rational double point singularity of the Weierstrass model to
a worse singularity (Sections \ref{sec:geomsubs}, \ref{sec:bran},
\ref{sec:quat}).
 We show how such degenerations are naturally associated to certain
representations of the group $G$ (Sections \ref{sec:bran} and
\ref{sec:quat}).

Conversely, in Section \ref{sec:sadov} we show how the assigned
representations
and
re\-pre\-sent\-a\-tion-the\-o\-re\-tic facts  determine the geometry
of the Calabi--Yau and the degenerations which can occur.
We then derive the formula for {\boldmath  $\mathcal R$}, in terms of $G$
(which is associated to
the ``generic" singularity) and its representations (associated to
the ``non-generic" singularities).

 Only a limited set of representations occur, and only certain of the
 ``non-generic
singularities:'' many others are in fact excluded by the assumption that
$\pi: X \to B$ is a (flat) elliptic Calabi--Yau fibration.

 The computation of
{\boldmath  $\mathcal R$} is slightly different in the case of the simply
laced  exceptional groups (including those from Deligne's ``exceptional
series''):
here we obtain a novel relation between the Coxeter numbers
and the rank of these groups (Section \ref{sec:deligne}).

\medskip

{\bf Acknowledgments:} We would like to thank
P. Aspinwall,
M. Finkelberg, K. Intriligator,
S. Katz, M. Larsen, D. L\"ust, W. Nahm, B. Ovrut,
M. G. Rossetti, K. Uhlenbeck, and K. Wendland for useful discussions and
encouragement while this paper was in preparation.

\section{The group $G$} \label{sec:G}

\begin{definition}\label{1}
An {\em elliptic Calabi--Yau
threefold with section}\/ is a proper, flat map $\pi:X\to B$ from a
nonsingular
projective complex threefold $X$ with trivial canonical bundle to a
nonsingular
surface $B$, whose general fiber is an elliptic curve, and which admits a
section $\sigma:B\to X$.  (During certain parts of our discussion,
we shall also assume
 that   the rank of the Mordell--Weil group $MW(X/B)$ of the elliptic
fibration is
zero.)

\end{definition}

Any such $X$ is a resolution of a possibly singular, Weierstrass model
$\overline{\pi} : W \to B$
\cite{Na1,G}.
$W$ can be described (locally) by a ``Weierstrass
equation''
\begin{equation}\label{W} y ^2= x^3 + fx +g, \end{equation}
where $f$ and $g$ are sections of line bundles on the base $B$.

\begin{lemma}\cite{BIKMSV,AKM}
In this set up we can naturally associate a reductive
Lie group $G$ to the fibration
as follows.  Let $[E]^\perp$ be the orthogonal complement within $H_4(X)$
of the elliptic fiber $E$, and let $\Lambda$ be the cokernel of the natural
map
$$\pi^{-1}: H_2(B) \to [E]^\perp.$$
Then $\Lambda$ serves as the coroot lattice of $G$, and $\Lambda\otimes
U(1)$ serves as the Cartan subgroup.  Moreover, to each component of the
discriminant locus is associated a local factor of the group, determined by
the generic Kodaira fiber along that component and by the monodromy, as
follows:

\medskip

\begin{tabular}{l|ccccccc|}\hline
\multicolumn{2}{|l|}{generic Kodaira fiber} & $I_n$ & $I_n$ & $II$ &  $III$ &
   $IV$ & $IV$   \\ \hline
\multicolumn{2}{|l|}{monodromy}      & $\Z_2$ & $\{e\}$ & $\{e\}$ & $\{e\}$ &
   $\Z_2$ & $\{e\}$  \\ \hline
\multicolumn{2}{|l|}{local group factor} & $\sp([\frac n2])$ & $\su(n)$ &
   $\{e\}$ & $\su(2)$ &    $\sp(1)$ & $\su(3)$ \\ \hline
\multicolumn{7}{c}{\quad} \\ \hline
\quad\qquad\qquad
&$I_0^*$ & $I_n^*$ & $I_n^*$ & $IV^*$ & $IV^*$ & $III^*$ & $II^*$\\ \hline
& $\Z_3$ or $\mathfrak S_3$ & $\Z_2$ & $\{e\}$ & $\Z_2$ & $\{e\}$& $\{e\}$&
   $\{e\}$\\ \hline
&$G_2$ & $\so(2n{+}7)$ & $\so(2n{+}8)$ & $F_4$ & $E_6$ & $E_7$ & $E_8$\\ \hline
\end{tabular}

\end{lemma}

 \begin{proof}
To specify a connected reductive group, it is enough to specify a compact
torus together with the collection of characters of that torus which will
serve as the weights for the semisimple part of the group.  The torus in
turn can be described as $\Lambda\otimes U(1)$ for some lattice $\Lambda$,
which is the form used in the statement of the lemma.  (This choice of
torus is dictated by physical considerations.)

To complete the specification, the weight spaces must be given.
 The Weierstrass model is singular along a (reducible) curve
$C$; the general singularity over each irreducible component of $C$
is a rational double point \cite{kod}. Let us consider the intersection
configuration
of the exceptional curves and the exceptional divisors on $X$.
 In most cases, the intersection matrix is (up to a sign) the unique Cartan
matrix of a Lie algebra $\mathfrak g$;
here we also find the non-simply laced algebras, as the exceptional curves
might
undergo a monodromy transformation as they move in the exceptional
divisors along the curve $C$. In some cases a more delicate argument is
needed \cite{AKM}.

\end{proof}

Note that the group is semisimple precisely when $MW(X/B)$ has rank $0$.

\begin{corollary}\label{def:g'} Let  $G$ be a non-simply laced group
(a local factor of the entire gauge group)
associated to the singularities over a curve $C$ in the discriminant,
as in the above proof. Then the exceptional curves in one homology
class are parameterized by a curve $C'$, a finite branched cover of $C$.
 The cover is of degree $2$ unless $G$ is locally isomorphic to $G_2$;
in the latter case, the degree of the cover is $3$.
\end{corollary}
\begin{definition} With the notation of the above corollary, we write
$g(C')=g'$.
\end{definition}
\begin{remark} If $B$ is ruled, from $X$ and the group $G$
we can construct a {K3} surface $S$ with a gauge bundle $H$,
the ``heterotic dual'' of $X$. Many of the physics predictions stated in
this paper
were originally derived by analyzing this  duality.
\end{remark}

\begin{definition}\label{def:2} Let $\bsa \subset B$ be the  ramification
locus of $\pi$.
 $\bsa$ is a divisor. We write
\begin{equation} \label{ram}
 \bsa = \bso \cup _{i \geq i} \bsi,
\end{equation}
 where $\bso$ is
the (possibly reducible) component over which the ``general'' singular
fiber is a node (Kodaira type $I_1$), and each $\bsi$ is an irreducible
component of $\bsa \setminus \bso$.
\end{definition}

\begin{remark}\label{CY} Since  $X \to B$ is an elliptic Calabi--Yau
\cite{Na1}, then:
$$\mathcal O _X \sim K_X = \pi ^*( \kb + \frac{1}{12}\bsa), \text{ and }
  \bsa \in |-12K_B|.$$
\end{remark}
 The choice of the notation and of the indices in Definition \ref{def:2}
is motivated by the gauge
group associated to the resolution of the general singular point of $\bsi$:
we denote in fact by $G_i$ this group and because $G_0$ is trivial, the
relevant groups
 are $G_i, \ i \geq 1$.

\section{A first look at {\boldmath  $\mathcal R$}: when {\boldmath
$X=W$}.}\label{sec:firstlook}

Let us define the fundamental invariant:
\begin{definition}\label{R}
{\r}$= \f \c(X) + 30 K_B^2.$
\end{definition}
 The following holds:

\begin{theorem}\label{th:svw} Under the hypothesis of Definition \ref{1},
suppose that in addition
$X=W$ is a smooth Weierstrass model. Then {\r}=0.\end{theorem}
\begin{proof}
Following  the algorithm
provided in \cite{SVW} one can in fact show that
$$\c (X) = - 60 \kb ^2.$$
This statement is also buried in the proof of  \cite[(3.10)]{Na1}; even though
the author claims it only for $B = \p2, \ \mathbb F_1, \bpt$.
\end{proof}

We should also point out  that {\r=0} is not a sufficient condition
for $W$ to be smooth (see also Section \ref{sec:isolated}).
 We will give an alternative proof of this theorem in Corollary \ref{svw2};
this proof will follow the mathematical ideas that arise
from considering the ``vanishing of the anomalies''
in string theory. The following definition is also motivated by the
physics literature:
\begin{definition} \label{anom1}
$ H_{ch}=$ {\r} $+ \dim (G)- \rk (G).$\end{definition}
Theorem \ref{th:main} describes explicitly $H_{ch}$.
 We discuss in detail the motivation behind this definition
and a dictionary between the geometry and the quantum field
theory in \cite{AD}.

\section{A look from string theory: \\ Gauge theory on $X$ and vanishing of
the anomalies} \label{sec:physics}

Before studying other properties of the invariant {\r} defined in
the previous section, we look at the gauge theory interpretation
of our setup.
{}From this point of view we consider
a gauge theory on $X$ (coupled to gravity), in which certain coefficients
of the curvature
are required to vanish: these are the ``anomaly cancellations'' whose
geometric
counterparts are a formula for {\r} and certain geometric constraints.
We will see that a formula for {\r} is the geometric
counterpart
of the first anomaly cancellation (Theorem \ref{th:s1} and \cite{AD}); the
second anomaly
cancellation and the corresponding geometric constraints
will be discussed in Section \ref{sec:sadov}.

\subsection{From elliptic fibrations to gauge theory}

 When one of the ``type II string
theories'' is
formulated on a ten-manifold
of the form $M^{3,1}\times X$ with $X$ a Calabi--Yau
threefold
and $M^{3,1}$ a flat spacetime of dimension four, the resulting theory has
a low energy
approximation which takes the form of a
four-dimensional quantum field theory with quite realistic physical
properties
(depending on certain properties of the Calabi--Yau threefold).

Elliptic Calabi--Yau
threefolds with a section $\pi:X\to B$ have also
been used in a different way in string theory.
We can ask what happens to the type IIA theory
in the limit when the
Calabi--Yau metric on $X$ is varied so that the fibers of the map $\pi$
shrink to zero area.  It turns out that the resulting physical theory has
a
low energy approximation which takes the form of a {\it six-dimensional}\/
quantum field theory.  This limiting theory can also be described more
directly, in terms of the periods $\tau(b)$ of the elliptic curves
$\pi^{-1}(b)$, regarded as a multi-valued function on $B$.  The type IIB
string theory is compactified on $B$ with the aid of this function, using
what are known as D-branes along the discriminant locus of the map $\pi$.
(This latter approach is known as ``F-theory.'')

{\bf Fact:}
This six-dimensional quantum field theory includes gravity as well as
 a gauge field theory whose gauge group
is the group $G$  defined
in Section \ref{sec:G}; in
order to be a consistent quantum theory, the ``anomalies'' of this
theory must vanish.

\subsection{Curvatures, anomaly polynomial and traces}

Schwarz shows \cite{Schwarz} that these models ($N{=}1$ theories in
six dimensions
with a semisimple group  $G$) are constrained by anomaly cancellation.
The anomaly is characterized by an eight-form made from curvatures
of the Levi--Civita connection and of the gauge connection.
This eight-form is naturally defined on an auxiliary eight-manifold
$Y$.\footnote{$Y$ is a manifold with boundary, whose boundary is the
product of $S^1$ and the original six-dimensional spacetime.}

If we have a manifold $Y$ equipped  with a principal $G$-bundle
$\mathcal{G}$ (the
``gauge bundle''), then the curvature $F$ of the gauge connection
is an $\ad(\mathcal{G})$-valued
two-form, where each fiber $\ad(\mathcal{G})_x$ of $\ad(\mathcal{G})$ is
isomorphic to the Lie
algebra $\mathfrak{g}$ of $G$, with $\mathcal{G}_x$ acting on
$\ad(\mathcal{G})_x$ via the adjoint action of $G$ on $\mathfrak{g}$.
Similarly, if $Y$ is equipped
with a (pseudo-)Riemannian metric, then the curvature $R$ of the
Levi--Civita connection is a two-form taking values in the endomorphisms
of the tangent bundle.

The {\bf ``anomaly polynomial''} is a differential form on $Y$
which involves
expressions like $\tr R^k$ and $\Tr_\rho F^{k}$, where $\rho$ is some
representation of
the Lie algebra.  These expressions are to be interpreted as follows:  the
representation $\rho$ can be regarded as a homomorphism $\rho:\mathfrak{g}\to
\End(V)$ for some (complex) vector space $V$.  As an endomorphism of $V$,
$\rho(F_x)$
can be raised to the $k^{\text{th}}$ power; the resulting endomorphism
$\rho(F_x)^k$ of $V$ has a {\it trace}, which we denote as
$$\Tr_\rho F^k=\trace_V \rho(F)^k.$$
Although this expression might have depended on the choice of isomorphism
to $\mathfrak{g}$, in fact it is invariant under the adjoint action of $G$
on $\mathfrak{g}$ and so is independent of choices.

Similarly, the expressions $\tr R^k$ are evaluated with the help of the
``vector'' representation of the corresponding orthogonal group.

\subsection{Vanishing of the anomalies}

The first requirement is the vanishing of the coefficient of a certain
curvature term,
which imposes restrictions on the choice of the group  and its ``matter''
representations
which can occur.
In \cite{AD} we discuss extensively the geometric realization of
this formula, when F-theory is compactified on an elliptic Calabi--Yau
threefold:
\begin{theorem}\label{th:s1}\cite{Schwarz, MVII, AD}
The anomalies are characterized by an eight-form, made from curvatures
and gauge field two-forms.
One requirement is the vanishing of the coefficient of the curvature term
$\tr R^4$, where $R$ is the curvature of the Levi--Civita connection.
In our geometric set up this  leads to
\begin{equation}
\text{{\r}}=H_{ch}- \dim (G) + \rk (G),
\end{equation}
where $H_{ch}$ involves the dimension of certain representations of the
group $G$.
\end{theorem}
  The representations which occur in $H_{ch}$ are well defined in terms of
quantum field theory;
this motivates our Definition \ref{anom1}. Theorem \ref{th:main}
identifies these representations
in our geometric set up.
 We will analyze the geometric counterpart of the
following statement (``The generalized Green--Schwarz mechanism'')
in Section \ref{sec:sadov}:

\begin{theorem}\label{th:s2}\cite{GS,Sagnotti,Schwarz} Let us assume that
$G$ is semisimple, and so in particular that $\rk{MW}=0$ and that
$G$ is
locally isomorphic to $\prod_i G_i$, where $G_i$ are
simple groups.
If the requirement specified in
 Theorem \ref{th:s1} holds then the remaining terms of the anomaly
polynomial
(in a suitable normalization \cite{Sadov}) are:
\begin{equation}\label{I} \frac{9-n_T}8 (\tr R^2)^2 + \frac{1}{6} \tr R^2 \sum
X_i ^{(2)} -\frac{2}{3} \sum X_{i} ^{(4)} + 4\sum _{i < j
 }  Y_{ij}
\end{equation}
where $n_T$ denotes the number of ``tensor multiplets'' (which coincides with
$h^{1,1}(B)-1$ in our theories), and where
\begin{align*}
X_i ^{(n)} &= \Tr_{\adj} F_i ^n - \sum_\rho {n_\rho} \Tr _\rho F^n _i\\
Y_{ij}&= \sum _{\rho,\sigma} {n_{\rho\sigma}} \ \Tr _{\rho} F^2 _i \  \Tr
_{\sigma} F^2 _j.
\end{align*}
$\Tr_{adj}$ means the trace in the adjoint representation, $\Tr _\rho  $
denotes the trace in the representation $\rho$
of the simple group $G_i$ (see the above subsection),
$n_\rho$ is the multiplicity of the representation $\rho$ of $G_i$ in the matter
representation,\footnote{In the physics literature one says that there are
``$n_\rho$ hypermultiplets in the representation $\rho$.}
and $n_{i,j}$ is the 
multiplicity of
the representation
$(\rho,\sigma)$ of $G_i \times G_j$.

The Green--Schwarz cancellation mechanism (in the generalized form due to
Sagnotti \cite{Sagnotti}, see also Sadov \cite{Sadov})
says that the anomalies can be cancelled provided that \eqref{I} can be
written in the form:
\begin{equation}
\label{eq:factored}
\left({\bf s} \tr R^2 + \sum {\bf t_i} \tr F^2_ i \right)\cdot
\left( {\bf u} \tr R^2 +
\sum {\bf v_i} \tr F^2 _
i \right),
\end{equation}
where ${\bf s}$, ${\bf t_i}$, ${\bf u}$, and ${\bf v_i}$ are divisors on
the base $B$ (which
correspond to ``tensor multiplets'' in the physical theory), the product is
calculated using the intersection pairing on $B$, and $\tr F^2
_i$ is evaluated in an appropriate ``fundamental''
representation of $G_i$.
\end{theorem}

In the case that $G$ is not semisimple, the anomaly polynomial is also
known, but it is much more complicated.  In this paper, we will only
consider the anomalies associated to the Green--Schwarz mechanism in the
semisimple case.

\section{ About {\boldmath $H_{ch}$}: a look from the physics
literature.}\label{sec:charged}

We state some of the physics predictions on $H_{ch}$, based on Schwarz's
analysis; these predictions motivated our geometric definition of $H_{ch}$
(see Definition \ref{anom1}).  (We have only described a small number of
the predictions which appear in the physics literature---others can be
found in \cite{witten:MF,Asp,Sadov,KV,IMS,CPR,DE}.\footnote{Note in particular
that \cite{CPR} used anomaly cancellation---as we do---to 
make and verify predictions
about $H_{ch}$.})

\medskip

{\bf Case 0:}
If $W=X$ is a smooth Weierstrass model, that is $G=\{e\}$ and $\dim (G)-
\rk (G)=0$,
then the quantum field theory tells us that $H_{ch}=0$ and {\r} $=0$, as
$H_{ch}$ is the (sum of) dimensions of certain irreducible representations
of $G$ (see \ref{th:s1}).
 This is in agreement with Theorem \ref{th:svw}.

\medskip

{\bf Case I:}
If  $W$ is singular along a single, smooth curve of genus $g$ of
$A_{N-1}$ singularities everywhere, we know from Section \ref{sec:G} that
$G= \su (N)$.
The authors of \cite {KMP} show that under these hypothesis $$H_{ch}=
g(\dim (G)-\rk (G)),$$
and also state that the same should hold for any isolated curve.
In this case one would have: $$\text{\r} =(g-1)(\dim (G) - \rk (G)).$$

\medskip

{\bf Case II:} If the group is non-simply laced (see Section \ref{sec:G})
and $W$ is singular along a unique curve $C$ of genus $g$, then
some of the exceptional divisors in $X$
mapping to $C$ are ruled surfaces over a curve $C'$ of genus $g'$
(see Corollary \ref{def:g'}).  Assume that there are $B_1$ branch points
of the map $C'\to C$, and that all degenerations of the generic singular
fiber
occur at these branch points.
The authors of \cite{AKM} show that in most such cases
$$H_{ch}=g (\dim (G)-\rk (G)) + (g'-g)  \mbox{\rloc0},$$
where {\rloc0} is a constant (which corresponds to the ``charged
dimension'' of a certain representation $\rho_0$
of $G$---see Definition \ref{def:chargedim}).
In this case one would have:
\begin{equation}\label{eq:cdone}
\text{\r} =(g-1)(\dim (G) - \rk (G))+ (g'-g)\mbox{\rloc0} .
\end{equation}

In the case of  $I_{2k+1}$ with monodromy (yielding gauge group
$G=\sp(k)$),
this formula is modified to one which
involves $B_1$ as well:
\begin{equation}\label{eq:cdthree}
\text{\r} =(g-1)(\dim (G) - \rk (G))+ (g'-g)\mbox{\rloc0} +\frac12B_1(2k).
\end{equation}

\medskip

{\bf Case III:} If $W$ is singular along a single, smooth curve of genus $g$,
the singularities are generically of type
$A_{N-1}$ singularities, but they become of type $A_N$ at $B_2$ isolated
points: we know from Section \ref{sec:G} that $G= \su (N)$. The authors of
\cite{BIKMSV} and \cite{BKKM}
 show that under these hypothesis $$H_{ch}= g(\dim (G)-\rk (G)) + B_2 N.$$
In this case one would have: $$\text{\r} =(g-1)(\dim (G) - \rk (G)) + B_2 N.$$

\medskip

In Section \ref{sec:inter} we will prove that all of these predictions
hold and give a global explanation for the above formulas;
we will also derive the value of {\rloc0} (which depends on $G$).

\section{Working assumptions and (most of the) notation} \label{sec:wa}

Our basic strategy for verifying the formula for {\r} is as follows.
On the one hand, the Euler characteristic of $X$ can be calculated
exploiting the elliptic fibration, studying the various types of singular
fibers which can occur, and assigning to each a ``contribution'' to the
Euler characteristic.  First, the generic fibers make no contributions.
Second,
the fibers over the curves $\bsi$ make contributions which can be accounted
for in terms of the genus of $\bsi$ and of its monodromy cover as well as
the type of the Kodaira fiber.  This leaves the contributions from
intersection points of the $\bsi$'s, or from special points along the $\bsi$'s
at which the fiber becomes worse.

On the other hand, a parallel decomposition can be made of the
representation theory.  There are specific contributions to $H_{ch}$
which are associated to the various local
factors $G_i$ of the gauge group,
and depend on the genus of $\bsi$ and of its monodromy cover.  If
these are subtracted from our formula, what remains is a sum of
contributions from the intersection points of the $\bsi$'s, or from
special points along the $\bsi$'s at which the fiber becomes worse.

Thus, once the ``generic'' singularities have been matched up, the
verification can be reduced to a local question---for each type of
singular fiber, verify that its contribution to the Euler characteristic
is compatible with the assignment of a factor in the representation
to the fiber.

We will carry this out under some assumptions about the degenerations.
To simplify matters and isolate the core of the problem, we will
consider the case of a single non-abelian factor $G_1$ in the gauge
group.  We will also make some simplifying assumptions about which
degenerate fibers are allowed.  (The cases we consider can be extended
to a more general set-up: see  Remarks~\ref{rem:611}
and \ref{rem:87} as well as \cite{AD}.)
Our specific assumptions are
as follows (see Equation \eqref{ram}):

\begin{itemize}
\item[$(\bullet)$]
The locus of  enhanced gauge symmetry is over a unique
smooth curve $\bs1$.

\item[$(\bullet)$] The coefficients in the Weierstrass equation are
otherwise general;
following \cite{BIKMSV} we assume that the local equations can be
determined by the data in Table~\ref{table:1}.
\end{itemize}

\begin{proposition}\label{pro:sigma1}
Under these hypothesis the group $G$ alone determines the multiplicity $m$
of $\bs1 $ in $\bsa$
(see \cite[Table 2]{BIKMSV} and the Tables in Appendix I):
\begin{equation*} \bso + m \bs1 = \bsa \in |-12 K_B|.
\end{equation*}
Equivalently

\begin{equation*} \bsa \text{ is defined by the equation } s^m \sigma_0
,\end{equation*}
 where $ \sigma_0 $ defines $ \bso$ and $ \bs1$ is defined by $s=0.$
\end{proposition}

\begin{definition}\label{mufg}
We denote by  {\boldmath ${\mu (f)}$} and {\boldmath $ \mu (g)$} the
multiplicity of $f$ and $g$ resp.\
along $\bs1$, and by
 { \boldmath $\mu_P (f,g)$ }  the intersection multiplicity of
$f/s^{\mu(f)}$ and $g/s^{\mu(g)}$
at a point $P \in \bs1$.
\end{definition}

\begin{definition}  We denote by $X_{\bs1}$
 the singular fiber of
Kodaira type over the {\it general point of}\/ $\bs1$.

We denote by $\{ Q_1, \cdots Q_C \}$, the singularities of
$\bso $ {\it away}\/ from $\bs1$: these are cusps:
 $C$ is then the {\it number}\/ of cusps of $\bso$.

 We denote by $X_{\bso}$ the
 singular (nodal) fiber over the {\it general point}\/ of $\bso$ while
$X_C$ is the singular {\it (cuspidal)}\/ fiber over each point $Q_j$.

\end{definition}

 If $\bso$ and $\bs1$ are disjoint, all the degenerate elliptic fibers
are the ones described above; if $\bso \cap \bs1 \neq \emptyset$
there are other degenerate elliptic fibers, not necessarily of Kodaira
type, over each
 intersection point. A complete classification of such degenerations is not
available, except in the case of simple normal crossings \cite{Mi,Na2},
and the list of possibilities could be quite complicated.

 These points (the $P^i_\ell$ below) are exactly the singularities of
$\bsa$ along $\bs1$; the roots of $\sigma _0  \bmod s$ determine the
intersection
of
$\bs1$ and $\bso$.

In particular:
\begin{proposition}\label{link}\label{nu}
Our assumptions imply the following:
\begin{itemize}
\item[($\bullet$)]  The equation defining $\bso   \bmod s =0$
splits in the product of at most two factors:
$\beta ^{r_1} _1 \cdot \beta ^{r_2} _2 \equiv \sigma _1  \bmod s$.
Each $\beta _i$ is irreducible, and together with $r_i$ is
 determined by the choice of the group $G$  (see the Tables in Appendix I).
\item[($\bullet$)]  $\beta _i$ is smooth near $\bs1$: $(s, \beta _i)$ are
local coordinates around
each intersection point and $r_i$ is the intersection multiplicity of
$\beta _i$ with $\bs1$.
We write:
$$\bso \cap \bs1= \{ P_1 ^1, \cdots P^{B_1} _1, P^1 _2, \cdots P^{B_2} _2
\},$$ where $B_i$ is the {\it number}\/ of the distinct roots of $\beta _i$;
note that if $r_i=1$, the $P^{j}_i$'s are points of
simple normal crossings intersections.
\item[($\bullet$)] (Equivalently:) $ \bso \cdot \bs1 = (-12 K_B -m \bs1)
\cdot \bs1= r_1 B_1 + r_2B_2.$
\item[($\bullet$)] The degenerate elliptic fiber   over each point  $P^\ell_i$
and the local equation around $P^\ell _i$
does not depend on $\ell$, but only on $i=1,2$:
without loss of generality we write $X_{P_i}=\pi ^{-1} (P^i _ \ell)$.
We write the local equation in Table~\ref{table:2}.
\item[($\bullet$)] The intersection multiplicity $\mu_{P^\ell_i}(f,g)$ does
not depend on $\ell$, but only on $i=1,2$; we denote it by $\mu_i(f,g)$.
\end{itemize}
\end{proposition}
\begin{proof} It follows from \cite{BIKMSV}.
\end{proof}

\begin{proposition}\label{pro:data}

\begin{itemize}\item If $G= \so (k)$ ($k > 8$), and $m$
is as defined in Proposition \ref{pro:sigma1}  then  $\rk (G)+2=m$;
\item if $G= e, \su(2), \su (3), \so(8), E_6, E_7, E_8$ , and $J$ is regular,
then $\rk (G)+2=m$, $\mu_P (f,g)=0$;
\item
if $G= \su(n)$ and $J$ has a pole along $\bs1$, then $\rk (G)+1=m$, $\mu
(f), \ \mu (g)=0$,
\end{itemize}
\end{proposition}

Note that $G= \{e\}$ corresponds here to the Kodaira fiber of type $II$.

\begin{proof} It can be verified by inspection and explicit computations.
\end{proof}

We list the values of $m, r_1$, and $r_2$  in the Tables in Appendix I.

\section{De constructing {\boldmath  $\mathcal R$}} \label{sec:ec}

In this section we set up an algorithm to compute  {\boldmath  $\mathcal
R$}, the fundamental
invariant defined in \ref{R}.

We break up the contributions to $\c$ as follows:
\begin{lemma}\label{c1} The following lines  add to the topological Euler
characteristic of $X$:
\begin{align*}\label{em} \f \c (\cup _{i,\ell} \  \pi ^{-1}(P_\ell ^i))=&
 \f \c (X_{P_1}) \cdot B_1 + \f \c (X_{P_2} ) \cdot B_2 \\
  \f \c ( \pi ^{-1}({\bs1 \setminus \cup _{\ell,i} P_{\ell} ^i}))=&
 \c (X_{\bs1} )\cdot  [1- g(\bs1)- \f B_1 -\f B_2] \\
 \f \c ( \pi ^{-1} (\bso \setminus \cup _j Q _j
\setminus \cup _{\ell,i} P_{\ell} ^i)=&
\f  \c (\bso \setminus \cup _j Q _j
\setminus \cup _{\ell,i} P_{\ell} ^i)\\
 \f \c (\cup _j  \pi ^{-1} (Q_j))=& \f \c (X_C) \cdot C
\end{align*}
\end{lemma}
\begin{proof}
We compute the Euler characteristic of $X$ via the structure of elliptic
fibration
(the Euler characteristic of the general fiber is zero) and Mayer--Vietoris'
sequence.
\end{proof}
 Now we want  to effectively calculate each contribution in
the above equations in terms of quantities which depend on the
singularities along $C$ and the group $G$.

Note also that the singularities along
$C$ are determined by the geometry of the discriminant locus
on $B$; this is in turn determined by the intersections of a section of
some multiple
of  $K_B$ with $\bs1$ (see Remark \ref{CY}).

\subsection{\bf De constructing $\ \c (\bso \setminus \cup _j Q _j
\setminus \cup _{\ell,i} P_{\ell} ^i)$:}\label{subsec:decec}

\noindent\phantom{XX}

We start with the following definitions:
\begin{definition}\label{alfa}{\bf (Defining {\boldmath $\alpha _j$}.)}
If $\phi _1 : B_1 \to B$ is the blow up of a point $P \in D$, with
exceptional divisor $E$ and
$$D_1 = \phi _1 ^*(D) - \a _1 (P) E, \text{ its strict transform.}$$
\end{definition}

\begin{corollary} With the above notation:
$$(K_{B_1} + D_1) \cdot D_1 = (K_{B} + D) \cdot D - \a _1(P)\cdot (\a _1
(P)-1).$$
\end{corollary}
In particular, if $P$ is a smooth point of $D$, $\a _1(P)=1$; if $P$ is a
cuspidal point of $D$, $\a _1(P)=2$.

\begin{definition}\label{fi}Let $\phi ^i_{n(i)} \cdots \phi^ i _v \cdots
\phi ^i_1$ be the embedded
resolution of $\bso$ around the point $P ^i_1 $ and $ \{ \a^1 _v \} := \{
\a _v (P_1 ^i) \}_{v=1} ^{n (i)}$
the collection of the integers as in \ref{alfa}
($v$ depends on $i$, but we believe the distinction is clear.)
Let us define:
 $$\epsilon_j=[\sum \alpha ^j_v (\alpha ^j _v -1) - \# \phi ^{-1}(P^k_j)] .$$

\end{definition}

 If $\ P^i_1$ is a smooth point,  $\epsilon _1=-1$;
for the cuspidal points $Q_j$ we have  $\epsilon=1$.

\begin{corollary}
$$\c (\bso \setminus P_j) =
\c (\tilde {\bso}) \setminus \#  \phi ^{-1}( P_j)= - (K_{B} + \bso) \cdot
\bso +
\epsilon_j $$
\end{corollary}
\noindent

\begin{proposition}\label{cusps} With the above notation we have:
\begin{equation*}
\begin{array}{rlc}
\c (\bso  \setminus  \cup _j Q _j
\setminus \cup _{\ell,i} P_{\ell} ^i) +
11 \cdot 12 \kb ^2  =&m \kb \cdot \bs1 + 2m \bs1 \cdot \bso +
m^2 \bs1 ^2 +\\
&+ B_1\epsilon _1 + B_2 \epsilon _2 +C.
\end{array}
\end{equation*}
 $\epsilon _1$ and $\epsilon _2$ are defined
in \ref{fi}; they are determined by non-generic
the singularities along $C$.
\end{proposition}

\begin{proof}
{}From the previous corollary we have:
$$ \c (\bso \setminus \cup _j Q _j
\setminus \cup _{\ell,i} P_{\ell} ^i) = -(\kb + \bso) \cdot \bso + \epsilon
_1 B_1
 + \epsilon _2 B_2 + 1 C.$$

Note that  $\bso \in |-12 \kb -m \bs1|$, which gives:
$$ -(\kb + \bso) \cdot \bso = -11 \cdot 12 \kb ^2 + m \kb \cdot \bs1 + 2m
\bs1 \cdot \bso
+ m^2  \bs1 ^2 . $$
 Substituting this in the above equation we obtain the statement of the
proposition.
\end{proof}

\subsection{\bf De constructing $C$, the number of cusps:}\label{subsec:decc}
\noindent \phantom{KEJWSHT}

\begin{lemma}
 $f $ and $g$ then have
$$(-4 \kb - \mu (f) \bs1)( -6 \kb - \mu (g) \bs1)=24 \kb ^2 + \{ 4 \mu (g) +
6 \mu (f)\} \ks + \mu (f) \mu (g) {\bs1}^2$$
intersection  points, counted with multiplicity ($f \in |-4 \kb|, \ g \in
|-6 \kb|$).
\end{lemma}

\begin{proposition}\label{c} Then number of cusps $C$, away from $\bs1$ is:
$$C= 24 \kb ^2 +  \{ 4 \mu(g) + 6 \mu(f)  \} \ks
-\mu_1(f,g)B_1-\mu_2(f,g)B_2 +
\mu (f)\mu (g) {\bs1}^2,$$
 where $\mu_i (f,g)$   are defined in  Propositions
\ref{mufg} and
\ref{nu}.
\end{proposition}

 $\mu_i (f,g), \ \mu (f)$, and $\mu (g)$  depend on the equation
\eqref{Wg} and are determined by the (non-generic and generic)
singularities along $C$ .
\begin{proof}
$C$ is the number of cusps away from $\bs1$; our assumptions in Section
\ref{sec:wa}
imply that the cusps are determined by the common zeroes of the polynomials
$\{f=g=0 \}$  away from $\bs1$ (these are ordinary vanishing, see equation
\eqref{fg}).
  $f$ and $g$ might also vanish along $\bs1$, of orders
$\mu (f)$ and $ \mu (g)$; $f \bmod  s$ and $g \bmod  s $ might have a common
zero along $\bs1$. The multiplicities of these latter zeros are measured by
$\mu_i(f,g)$. (See Appendix I.)
\end{proof}

\begin{proposition}\label{rev} Using the formulas \eqref{cusps} and
\eqref{c} derived above, we re-arrange the contribution to $\c(X)$ in
 \ref{c1} as follows:

$$\aligned
 \f \c (\cup _{i,\ell} \  \pi ^{-1}(P_\ell ^i)) \phantom{quagliuzzocanav}=&
 \f \c (X_{P_1}) \cdot B_1 + \f \c (X_{P_2} ) \cdot B_2 \\
\f \c ( \pi ^{-1}({\bs1 \setminus \cup _{\ell,i} P_{\ell} ^i}))
 \phantom{quagliuzzocan}=&
   ( g(\bs1)-1) (-m)- \f m B_1 -\f m B_2 \\
 \f \c ( \pi ^{-1} (\bso \setminus \cup _j Q _j
\setminus \cup _{\ell,i} P_{\ell} ^i)+ 54 \kb ^2 =
& \f  m \kb \cdot \bs1 + \f m^2 \bs1 ^2 + m \bs1 \cdot \bso \\
&+ \f \epsilon _1 B_1 + \f \epsilon _2 B_2 \\
 &+  \f [\{ 4 \mu(g) + 6 \mu(f)  \} \ks ] - \frac12\mu_1 (f,g) B_1\\
&- \frac12\mu_2 (f,g) B_2 + \f [\mu (f) \mu (g) {\bs1} ^2] \\
 \f \c (\cup _j  \pi ^{-1} (Q_j) \phantom{quagliuzzo}-24 \kb ^2 =& \{ 4
\mu(g) + 6 \mu(f)  \} \ks - \mu_1(f,g) B_1\\
&- \mu_2(f,g) B_2+ \mu (f) \mu (g) {\bs1}^2
 \endaligned$$
The entries in the left hand sides of the above equations add to {\boldmath
 $\mathcal R$};
 the coefficients on the right hand side
are determined by the singularities, generic and non-generic, along $C
\subset W$.
\end{proposition}

\begin{corollary}\label{svw2} If $X$ is a smooth Weierstrass model,
then {\r}=0.
\end{corollary}
One of the aims of this paper is to show that the  entries on the right hand
side are
a collection of dimensions of certain representations of $G$
(Main Theorem \ref{th:main})
which are determined by the singularities,
generic and
non-generic along $C \subset W$. In Section \ref{sec:sadov} we will
show
that also the converse is true, that is, the assigned representations
determine uniquely the geometry of $W$.

\begin{remark}\label{rem:611}
Note that the formula in Proposition~\ref{rev} admits an immediate
generalization
to cases in which there are more simple factors in the gauge group
(corresponding
to additional components of the discriminant).  A somewhat more involved
notation is required, to handle possibilities of singular curves $\bsj$
or intersections among several components, but the same geometric principles
we used above will lead to a formula of the same general type.
\end{remark}

\section{A second look at {\boldmath  $\mathcal R$}: $\bs1 \cdot \bso=0$
($\bs1$ is isolated).}\label{sec:isolated}

We have considered in Section \ref{sec:firstlook} the case
$X=W$, we consider now the case when $\bs1$ does not intersect
the rest of the discriminant locus: equivalently, $W$ is singular along a
single curve $C$ and the singularities are uniform along $C$.
This case was also considered in the physics literature, see Section
\ref{sec:charged}.

Here the computations are
simpler, and we can see clearly how by using the geometry of the
base we can write {\boldmath  $\mathcal R$} (that is,
the equation in \eqref{rev})
as a function of the singular locus and certain representations
of the group $G$. The first implication of the hypothesis is that
$J: B \dasharrow  {\mathbb P^1} $ is well defined around $\bs1$.
 By analyzing the vanishing of the anomaly we find a curious relation
between the Coxeter number and rank in the case of the ``exceptional
series'' of Deligne.

\begin{theorem}\label{sec:deligne} If  $\bs1$ does not intersect the other
components
of the discriminant locus then
$$\text{ \r } = (\dim (G) - \rk (G)) (g-1).$$
Note that $\dim (G) - \rk (G)= \dim \adj _G - \dim Ker (\adj _G)$.
\end{theorem}

\begin{proof}

\medskip

{\bf Case I: $J $ is regular along $\bs1$}
(simply laced groups in Deligne's  exceptional series.)

In this case $G= \{e\},\  \su(2), \su (3), \so (8), E_6, E_7, E_8$, which
(except for the trivial case) are precisely the simply laced groups in
Deligne's exceptional series.
 Here the singular fibers are of
types $II, III, IV, I^* _0, IV^*, III^*, II ^*$,
and $m= \c (X _{\bs1}) <12$.

\medskip

By assumptions,
 $B_1=B_2=0$ and  {\boldmath ${\mathcal R}$} becomes (\eqref{rev}):

$$\aligned
\f \c ( \pi ^{-1}({\bs1 \setminus \cup _{\ell,i} P_{\ell} ^i}))
 \phantom{quagliuzzocan}=&
   ( g(\bs1)-1) (-m) \\
 \f \c ( \pi ^{-1} (\bso \setminus \cup _j Q _j
\setminus \cup _{\ell,i} P_{\ell} ^i)+ 54 \kb ^2 =
&\f m \kb \cdot \bs1 +  \f m^2 \bs1^2 \\
&+ \f \{ 4 \mu(g) + 6 \mu(f)  \} \ks \\&+ \f [ \mu(f) \mu (g) \bs1 ^2]
\\
 \f \c (\cup _j  \pi ^{-1} (Q_j) \phantom{quagliuzzo}-24 \kb ^2 =& \{ 4
\mu(g) + 6 \mu(f)  \} \ks
+   \mu(f) \mu (g) \bs1 ^2
 \endaligned$$

Now we use the geometry of the singularities of $W$:
$$ \aligned \kb \cdot \bs1 + \bs1^2 &= 2g(\bs1) -2 \\
  \bso \cdot \bs1=0 \ \leftrightarrow \ \-12\kb \cdot \bs1 &= m \bs1
^2\endaligned$$
By solving the system we have ($m < 12$):
$$ \aligned -\kb \cdot \bs1 &=  \frac{m}{12-m} 2(g-1) \\
 \bs1 ^2 &=  \frac {12}{12-m} 2(g-1) \endaligned$$

By substituting the above equation in the right hand side  of {\boldmath
$\mathcal R$},
we see that every term is a multiple of $g-1$.

Then:
$$\text{ \r }= \frac{6}{12 -m} \{ -2m + [m-2 \mu(g)][m-3 \mu(f)] + m^2 \}
(g-1)= \frac{6}{12 -m}m (m-2) (g-1),$$
in fact, by the definition of $m$, $[m-2 \mu(g)][m-3 \mu(f)]=0$ (see
Appendix I).

In this case we also have $\rk (G)= m-2$ (see Proposition \ref{pro:data}) and
thus:
$$\text{ \r } =\frac{6 (\rk(G)+2)}{10-\rk(G)}\rk(G)(g-1) .$$

Now, for these groups
$$\frac{6 (\rk(G)+2)}{10-\rk(G)}\rk(G)= h(G)= \rk(G) - \dim (G),$$
where $h(G)$ is the Coxeter number (see the following Lemma \ref{lemma:coxt}).

This is in agreement with the
expectations from physics (see Section \ref{sec:charged}) together
with Corollary~\ref{anom1}.

\medskip

{\bf Case II: $J $ has a pole along $\bs1$}
Since $J: B \dasharrow \pl$ is well defined around $\bs1$,
$ J_\infty \cdot \bs1=\kb \cdot \bs1=0$. This together with
the assumption
$0=\bso \cdot \bs1= (- 12\kb -m \bs1) \cdot \bs1 =-m \bs1^2$,
 implies
$2g-2=0$, that is, $g=1.$

 The substitution in  \eqref{cusps}
gives {\r} $=0$.
This is again consistent with the expectations in
Section \ref{sec:charged}.\end{proof}

\begin{lemma}\label{lemma:coxt} The Coxeter numbers of the simply laced
groups in  Deligne's exceptional series satisfy the relation:
$$ h(G)=\frac{6 (\rk (G)+2)}{10-\rk (G)}.$$
\end{lemma}

\begin{proof}
Case by case checking.
\end{proof}

 This adds to the numerology of the exceptional series
presented by Deligne in \cite{deligne}.

\section{Another look at \r : the Main Theorem.}\label{sec:inter}

In the discussion below,
we will describe the matter representation as a representation of
the Lie algebra $\mathfrak{g}$; it is in fact induced from
a representation of
the full gauge group $G$ associated to $X$.

\begin{definition} \label{def:chargedim}
Let $\rho$ be a representation of a Lie algebra
$\mathfrak{g}$,
with Cartan subalgebra $\mathfrak{h}$.
The {\em charged dimension of $\rho$} is
$(\dim \rho)_{ch}= \dim (\rho)- \dim (ker \rho|_{\mathfrak{h}})$.
\end{definition}
For example, if $\rho$ is the adjoint representation then
$$(\dim \adj) _{ch}= \dim \mathfrak{g}- \dim \mathfrak{h}= \dim G - \rk G.$$

\begin{theorem}\label{th:main}
Notation as in Section \ref{sec:wa}.
Then:

$$\text{\r} = (g-1) \dim (adj) _{ch} + (g' -g)  \dim
(\rho_0)_{ch}  + \sum _{P \in \mathcal A}
\delta_P
\dim(\rho_P)_{ch} , $$
where
 $\mathcal A = \{ P \in \bs1 \cap \bso $ such that the fiber over $P$
is of Kodaira type $\}$,
$g'$ is defined in \ref{def:g'},
the representations $\rho_P$ all come from a small list of representations
given in Table~\ref{tab:A},
and the coefficient $\delta_P$ is $\f$
if the representation is quaternionic
and is $1$ if the representation is real or complex.
(The quaternionic cases are labeled with $\f$ in the Table.)

In Table~\ref{tab:A}, we give the Kodaira type of the
general hyperplane section through the singular fibers
which occur under our hypotheses.  For each type of
singular fiber, we either list the associated representation
$\rho_j$, or (in the case of monodromy) we separate the
``non-isolated part'' of the representation and call it
$\rho_0$, listing any residual representation as $\rho_j$.
In addition, in a few
cases a representation occurs with multiplicity and (for later convenience
at the end of section~\ref{sec:sadov})
we identify an irreducible representation $\widehat{\rho}$ in the Table.

\end{theorem}

\begin{table}[ht]

\begin{center}

\begin{tabular}{||c|c||c|c||c||c|c||c||} \hline
Type & $G$ & $P_1$ & $P_2$ & $\rho_0$ &  $\rho_1$ & $\rho_2$ &
$\widehat{\rho}$  \\ \hline
$I_{1}$ & $\{e\} $ &$II$&$I_2$&&-- &\terminal& \\ \hline
$I_{2}$ & $\su (2) $ &$III$&$I_3$& &--&$\fund$ &\\ \hline
$I_{3}$ & $\su (3) $ &$IV$&$I_4$&  &-- &$\fund$& \\ \hline
$I_{2k} $, $k\ge2$ & $\sp (k) $&$I_{2k-4}^*$&$I_{2k+1}$&$\Lambda^2_0$&--&
    $\fund$  &\\ \hline
$I_{2k+1}$, $k\ge1$ & $\sp (k)$ &$I_{2k-2}^*$ &$I_{2k+2}$
    &$\Lambda^2+2\times\fund$ &$\frac12\fund$  &\terminal&$\fund$ \\ \hline
$I_{n}$, $n\ge4$ & $\su (n) $&$I_{n-4}^*$&$I_{n+1}$ && $\Lambda^2$&
   $\fund$  &\\ \hline
$II$ & $\{e\}$ &$III$&&  &\terminal &&  \\ \hline
$III$ & $\su (2)$ &$IV$&&&   $2\times\fund$& &$\fund$  \\ \hline
$IV$ & $\sp(1)$ &$I_0^*$&& $\Lambda^2+2\times\fund$&$\frac12\fund$ &
   &$\fund$ \\ \hline
$IV$ & $\su (3)$ &$I_0^*$&&& $3\times\fund$& & $\fund$ \\ \hline
$I_0^*$ & $G_2$ &$I_1^*$&&$\mathbf{7}$&-- & &\\ \hline
$I^*_{0}$ & $\spin (7) $&$I_1^*$&$I_1^*$&$\vect$&-- &$\spinrep$ & \\ \hline
$I^*_{0} $ & $\spin (8) $&$I_1^*$&$I_1^*$&&  $\vect$&$\spinrep_\pm$ &\\ \hline
$I^*_{1}$ & $\spin (9) $&$I_2^*$&$IV^*$&$\vect$&-- &$\spinrep$  & \\ \hline
$I^*_{1} $ & $\spin (10) $&$I_2^*$&$IV^*$&& $\vect$ & $\spinrep_\pm$  &\\
  \hline
$I^*_{2}$ & $\spin (11) $&$I_3^*$&$III^*$&$\vect$&-- &$\frac12\spinrep$  & \\
   \hline
$I^*_{2} $ & $\spin (12) $&$I_3^*$&$III^*$&&  $\vect$& $\frac12\spinrep_\pm
   $&\\ \hline
$I^*_{n}$, $n\ge3$ &$\so (2n+7)$&$I_{n+1}^*$&\nonmin&$\vect$&-- &\nonmin&\\
   \hline
$I^*_{n} $, $n\ge3$ & $\so (2n+8) $&$I_{n+1}^*$&\nonmin&&  $\vect$&\nonmin
   & \\ \hline
$IV ^*$ & $F_4$ &$III^*$&&$\mathbf{26}$&--& &\\ \hline
 $IV ^*$ & $E_6$ &$III^*$&&&  $\mathbf{27}$ & & \\ \hline
 $III ^*$ & $E_7$ &$II^*$&& &$\frac12\mathbf{56}$ &&\\ \hline
 $II ^*$ & $E_8$ &\nonmin&&&\nonmin  &  &    \\ \hline
\end{tabular}
\end{center}
\smallskip
\caption{The representations which occur under our ``generic''
hypotheses.}\label{tab:A}
\end{table}

\begin{remark} Our assumption of a smooth, flat elliptic fibration,
imposes restrictions on the type of degenerate singular fibers that might
occur:

\medskip

(i) If $\{e\}$ is associated to the Kodaira type fiber $II$,
there is a double point
singularity in the fiber over the simple normal crossings intersection
point of the two branches
($\bso$ and $\bs1$). This is terminal but not canonical,
leading to a smooth but {\it not flat} fibration and
a non-minimal Calabi--Yau threefold. We assume then that such points do not
occur:
the curve is isolated and the theorem holds
(see Theorem \ref{sec:deligne}).

\medskip

(ii) If $G=\{e\}$ (associated to the Kodaira type fiber $I_1$)
or
$G=Sp(k)$ (associated to the Kodaira type fiber $I_{2k+1}$),
the resolution of the generic singularities leaves a double point
singularity in the fiber over the simple normal crossings intersection
point of the two branches
($\bso$ and $\bs1$). In fact, if the equation is otherwise generic, then no
small resolution exists.
We assume here for simplicity that there are no such points.

\medskip

(iii) If $G$ is associated to the Kodaira type fiber $II^*$,
or $I^*_n, \ n \geq 12$, the equation of the Weierstrass model is
not minimal at the non-simple normal crossings intersection point of the
two branches
$\bso$ and $\bs1$. In order to resolve this singularity we would
need to blow up $B$ the basis of the fibration. In the resulting elliptic
fibration (still flat and Calabi--Yau), the two branches of
discriminant are separated.
 We assume then that such points do not occur.

\end{remark}

\begin{remark}\label{rem:notation}

We have used the following notation in Table~\ref{tab:A}:
\begin{itemize}
\item[($\bullet$)]
Cases with no small resolution are denoted ``\terminal'', and cases with
non-minimal Weierstrass model are denoted ``\nonmin''.
\item[($\bullet$)]
A dash denotes the trivial
representation, whereas a blank entry denotes a situation in which there
is no representation which belongs in that location.

\item[($\bullet$)]
The classical groups $\su(n)$, $\sp(n)$ have representations on
$\mathbb{C}^n$, $\mathbb{H}^n$, respectively, which are known as the {\em
fundamental representations}\/ and denoted by ``$\fund$''.  This representation
is quaternionic in the case of $\sp(n)$.
The second exterior power of the fundamental representation is denoted by
``$\Lambda^2$''. In the case of $\sp(n)$, $\Lambda^2$ is reducible and its
irreducible ``traceless'' part is denoted by ``$\Lambda_0^2$''.

\item[($\bullet$)]
The classical group $\so(n)$ has a representation on $\mathbf{R}^n$ called
the {\em vector representation}\/ and denoted by ``$\vect$''.  Its double
cover $\spin(n)$ has {\em spinor representations.} When $n$ is odd, there
is one spinor representation, of dimension $2^{(n-1)/2}$, denoted by
``$\spinrep$''.
When $n$ is even, there are two half-spinor representations, each of
dimension $2^{(n-2)/2}$, denoted by ``$\spinrep_+$'' and ``$\spinrep_-$''.
Note that the spinor or half-spinor representations are real if $n\equiv0, \pm1
\bmod8$, complex if $n\equiv\pm2\bmod8$, and quaternionic if
$n\equiv\pm1,4\bmod8$.

\item[($\bullet$)]
In the case of the exceptional groups, we label representations by their
dimension (given in boldface type).
\end{itemize}

\end{remark}

\begin{proof}
As we have already seen in Section \ref{sec:isolated},
the intersection numbers of the various parts of the discriminant in $B$
determine the geometry of $W$ and
 the choice of the
group $G$ and vice versa. Following Section \ref{sec:charged}, we write
all the terms in {\r} in Proposition \ref{rev}, as coefficients
of $g(C)$, the genus of the curve of singularities, the number of points
where the singularities are non-generic, and $g(C')=g'$,
when the groups are non-simply laced, and then interpret the results.
 The coefficients in \ref{rev} are determined by the group and
 the local geometry (the degeneration of the general
 rational double point) and are listed in Appendix I.
 We divide the proof in 3 steps.

$\bullet$ Step I (\ref{sec:geomsubs}):

We show how the geometry suggests the appropriate substitutions
for $\bso \cdot \bs1$, $\ks $,  $\cc$ and also $B_1$
if the group has monodromy branched at $B_1$ points.

 If $B_1=B_2=0$, then the substitutions are uniquely determined
(see Section \ref{sec:isolated}).

 In section \ref{sec:sadov} we show how these substitutions
 are equivalent to certain re\-pre\-sen\-ta\-tion-the\-o\-re\-tic facts.
 If $ G \neq \sp(k)$ or $\so(m)$, then after the substitutions
we obtain the data in  Table~\ref{tab:C}.  That is, the resulting
formula for
 {\r}\ can
 be written as a sum of local terms, associated to various points $P$,
 which can be collected into a formula of the form
\begin{equation}\label{eq:Rloc}
\mbox{\r} = (g-1)(\dim(G){-}\rk{G}) + (g'-g)\mbox{\rloc0}
+\sum_{j=1}^2 B_j \mbox{\rloc j}.
\end{equation}
The local contributions {\rloc j}
are recorded in Table~\ref{tab:C}.


\begin{table}[ht]

\begin{center}

\begin{tabular}{|c|c|c|c|c|c|} \hline
 Type & $ G $& $\dim{G}{-}\rk{G}$ &{\rloc0}&{\rloc1}&{\rloc2}\\ \hline
$I_{1}$ & $\{e\} $& $0$ &$0$&$0$  & \terminal  \\ \hline
$I_2$ & $\su (2)$ & $2$ &$0$&  $0$  & $2 $
  \\ \hline
$I_3$ & $\su (3)$ & $6$ &$0$&  $ 0$  & $3$    \\ \hline
$I_{2k}, k\ge2 $ & $\sp (k) $&$2k^2$  &$2k^2-2k$&$0$&$2k $   \\ \hline
$I_{2k+1}, k\ge1 $ & $\sp (k) $&$2k^2$  &$2k^2+2k$&$k$ & \terminal  \\
   \hline
$I_{n}, n\ge4$ & $\su (n) $& $n^2 -n$ &$0$& $\frac12(n^2-n) $  & $n $ \\
  \hline
$II$ & $\{e\}$ & $0$ &$0$& \terminal &      \\ \hline
$III$ & $\su (2)$ & $2$ &$0$&$4$  &     \\ \hline
$IV$ & $\sp (1)$ & $2$ &$4$ & $1$  &    \\ \hline
$IV$ & $\su (3)$ & $6$ &$0$&  $9 $  &    \\ \hline
$I_0^*$ & $G_2$ &  $12$ &$6$&$0$ &    \\ \hline
$I^*_{0}$ & $\spin (7) $&  $18$ &$6 $&$0$ & $8 $    \\
  \hline
$I^*_{0}$ & $\spin (8) $&  $24$ &$0$&$8$ & $8$    \\ \hline
$I^*_{1}$ & $\spin (9) $&  $32$ &$8$&$0$  &$16$    \\ \hline
$I^*_{1}$ & $\spin (10) $&  $40$ &$0$&$10$ & $16$   \\ \hline
$I^*_{2}$ & $\spin (11) $&  $50$ &$10$&$0$  &$16 $    \\ \hline
$I^*_{2}$ & $\spin (12) $&  $60$ &$0$& $12$& $16$    \\ \hline
$I^*_{n}, n\ge3$ & $\so (2n{+}7) $&  $2(n{+}3)^2$ &$2n{+}6$&$0$  & \nonmin   \\
  \hline
$I^*_{n}, n\ge3 $ & $\so (2n{+}8) $&  $2(n{+}3)(n{+}4)$ &$0$& $2n{+}8 $   &
  \nonmin   \\ \hline
$IV ^*$ & $F_4$ &  $48$ &$24$&$0$ &  \\ \hline
 $IV ^*$ & $E_6$ & $72$&$0$& $27 $  &    \\ \hline
$III ^*$ & $E_7$ &  $126$ &$0$& $28 $    &   \\
  \hline
  $II ^*$ & $E_8$ & $240$ &$0$&  \nonmin  &   \\
  \hline

\end{tabular}

\end{center}
\smallskip
\caption{The local contributions to the invariant {\boldmath
$\mathcal R$}.}\label{tab:C}

\end{table}

 In the cases $G = \sp(k)$, $G=\so(m)$, there are choices in making
 the substitutions but if a careful choice is made we can again
 write things in the form \eqref{eq:Rloc}
  (see also Section
 \ref{sec:sadov} for a better interpretation).
 
As we will point out in Remark~\ref{rem:87} below, the substitutions
can be formulated in a very general way which allows them to be applied
in cases beyond the specific ones considered here \cite{AD}.

$\bullet$ Step II (\ref{sec:bran}):
We show how we can naturally interpret the entries in Table~\ref{tab:C}
as charged dimensions of certain representations (multiplied by the
coefficient $\delta$), given in Table~\ref{tab:A}.  That is,
$\mbox{\rloc j} = \delta_j \dim(\rho_j)_{ch}$.
  If $p$ is not a branch point, then the (resolution of the)
  general elliptic surface
  through $P$ can be associated to a group $G'$ containing $G$, and
  the representation is obtained
   via the branching rules for the adjoint representation of $G'$.

  If $G$ is non-simply laced, then we consider $G \subset G'$, $G'$
  simply laced, and we use again the branching rules.
(This gives the representation-theoretic interpretation of the number
``{\rloc0}'' from equations
\eqref{eq:cdone},
\eqref{eq:cdthree}.)

  $\bullet$ Step III (\ref{sec:quat}).
   Finally
   we show how the
  number $\delta$ can be derived from the geometry of the degeneration
  of the general double point to the singularity over $p$.
\end{proof}

\subsection{Step I: The substitutions}\label{sec:geomsubs}

\begin{proposition}\label{le:su} Assume that resolution of the curve of
singularities
$C$ leads to a non-simply laced group $G$, as in \ref{def:g'}.
 Namely, some of the
exceptional divisors are ruled over a curve $C'$, which is
a finite cover of $C$ of degree $d=2,3$ ($3$ if and only if $G=G_2$),
ramified at $B_1$ points. Write $g'=g(C')$, then:
\begin{equation*}
B_1= 2(g'-g)-(2d-2)(g-1)
\end{equation*}
\end{proposition}
\begin{proof} The statement follows from Hurwitz's formula.
\end{proof}
\begin{proposition}\label{prop:85}
 Following the notation in Section \ref{sec:wa},
we have:
 $$\bs1 \cdot \bso = r _1 B_1 + r _2 B_2.$$
If the group $G$ is non-simply laced, then
$$ \bs1 \cdot \bso = r_2 B_2 + 2r_1 (g'-g) -(2d-2)r_1 (g-1),$$
if there are $B_1$ branch  points of: $C' \to C$.
\end{proposition}

\begin{proposition}\label{prop:86} The appropriate substitutions for
$-\ks$ and $\bs1 ^2$ are the ones given in Table~\ref{tab:B}.
\end{proposition}

\begin{table}[ht]

{\renewcommand{\arraystretch}{1.2}
\begin{center}
\begin{tabular}{|c|c|c|c|} \hline
Type & $ G $  & $-\ks $ & $ {\bs1 }^2$ \\ \hline
$I_{2} $ & $\su (2) $&  $\frac12B_1$ & $2(g-1)+\f  B_1$\\ \hline
$I_{2k}$, $k \geq 2$ & $\sp (k) $& $-(g-1)+(g'-g)$ &  $(g-1)+(g'-g)$ \\ \hline
$I_{2k+1}$, $k \geq 2$ & $\sp (k) $& $-(g-1)+(g'-g)$ &  $(g-1)+(g'-g)$ \\
  \hline
$I_{n}$, $n \geq
  3$ & $\su (n) $&  $ B_1$ & $2(g-1)+ B_1$\\ \hline
$II$&$\{e\}$&$\frac25(g-1)+\frac15B_1$&$\frac{12}5(g-1)+\frac15B_1$\\ \hline
$III$ & $\su (2)$&$\frac23(g-1)+\frac13B_1$& $\frac83(g-1)+\frac13B_1$
  \\ \hline
$IV$ & $\sp(1)$&$\frac12(g-1)+\frac12(g'-g)$ &$\frac52(g-1)+\frac12(g'-g)$
  \\ \hline
$IV$ & $\su(3)$&$(g-1)+\frac12B_1$ &$3(g-1)+\frac12B_1$  \\ \hline
$I_0^*$ & $G_2$  &$ \frac43(g-1)+\frac13(g'-g)$ &  $
  \frac{10}3(g-1)+\frac13(g'-g)$  \\ \hline
$I_0^*$ & $\spin(7)$ &$\frac53(g-1)+\frac13(g'-g)+\frac13B_2$&
$\frac{11}3(g-1)+\frac13(g'-g)+\frac13B_2$\\ \hline
$I_0^*$ & $\spin(8)$ &$2(g-1)+\frac13B_1+\frac13B_2$&
  $4(g-1)+\frac13B_1+\frac13B_2$\\ \hline
$I_n^*$, $n\ge1$ & $\so(2n+7)$ &$2(g-1)+B_2$& $4(g-1)+B_2$\\
 \hline
$I_n^*$, $n\ge1$ & $\so(2n+8)$ &$2(g-1)+B_2$& $4(g-1)+B_2$\\
 \hline
$IV ^*$ & $F_4$ & $3(g-1)+(g'-g)$& $5(g-1)+(g'-g)$ \\ \hline
$IV^*$ & $E_6$&$4(g-1)+B_1$ &$6(g-1)+B_1$  \\ \hline
$III^*$ & $E_7$&$6(g-1)+B_1$ & $8(g-1)+B_1$ \\ \hline
$II^*$ & $E_8$&$10(g-1)+B_1$ & $12(g-1)+B_1$ \\ \hline
\end{tabular}
\end{center}
}
\smallskip
\caption{The substitutions.}\label{tab:B}

\end{table}

\begin{proof}

\begin{itemize}
\item[(a)]
When $J$ is finite and there is no monodromy, i.e., cases $II$, $III$, $IV$,
$I_0^*$, $IV^*$, $III^*$, $II^*$ corresponding to the simply laced groups
 $\{e\}$, $\su(2)$,
$\su(3)$, $\so(8)$, $E_6$, $E_7$, $E_8$, in Deligne's exceptional series, then
the local geometry is given by the following equations:
$$ \aligned \kb \cdot \bs1 + \bs1^2 &= 2g(\bs1) -2 \\
 (-12\kb - m \bs1) \cdot \bs1 &= r_1 B_1+r_2B_2 \ ( \leftrightarrow \bs1 \cdot
\bso = r_1 B_1+r_2B_2), \endaligned$$
which can be solved since $m<12$:
$${
-\ks=\frac{2m(g-1)+r_1B_1+r_2B_2}{ 12 -m}  ;\quad
{\bs1}^2 = \frac{24(g-1)+r_1B_1+r_2B_2}{ 12 -m} .} $$
\item[(b)]
When $J$ is finite and there is monodromy, i.e., cases $IV$, $I_0^*$,
$I_0^*$, $IV^*$ corresponding to groups $\sp(1)$, $G_2$, $\so(7)$, $F_4$
which includes the remainder of Deligne's exceptional series,
the local geometry is the same but we also use Proposition \ref{le:su} to
eliminate $B_1$ in favor of $g'-g$:
$$ -\ks=\frac{(2m-2r_1(d{-}1))(g-1)+2r_1(g'-g)+r_2B_2}{ 12 -m}  ;$$
$$
{\bs1}^2 = \frac{(24-2r_1(d{-}1))(g-1)+2r_1(g'-g)+r_2B_2}{ 12 -m}.$$
\item[(c)]
If $G=\su(n)$, $n\ge3$, Table~\ref{table:1} in Appendix I tells us that
\begin{equation}\label{sub:1}
B_1=-\ks,
\end{equation}
where $B_1$ is the number of non-simple normal crossings
intersections.
The genus
formula then says that
\begin{equation}\label{sub:3}
\cc = 2(g-1) + B_1.
\end{equation}
(The case of $\su(2)$ is similar, using $B_1=-2\ks$.)
\item[(d)]
If $G=\sp([\frac n2])$, $n\ge3$, coming from $I_n$ with monodromy, then
$B_1=-2\ks$ so that
$$-\ks=(g'-g)-(g-1).$$
Combining this with the genus formula yields
$$\cc=(g'-1)+(g-1).$$
\item[(e)]
Finally, if $G=\so(2n+7)$ or $\so(2n+8)$ coming from $I_n^*$, $n\ge1$, then
$B_2=-2\ks-\cc$ which can be combined with the genus formula and solved to
give:
$$-\ks=2(g-1)+B_2; \quad \cc=4(g-1)+B_2.$$
\end{itemize}
\end{proof}

Step I now proceeds as follows: use the data in Tables~\ref{table:3} 
and \ref{table:4} in  Appendix I to evaluate
the ``local'' contributions to the Euler characteristic, in the formula for
{\r } given in Proposition \ref{rev}.
Then make the substitutions given in
Propositions
\ref{prop:85} and \ref{prop:86}
(supplementing them with Proposition \ref{le:su} if there is monodromy)
into the resulting formula; in all but a few cases (detailed below) this
yields a formula of the form
$$\mbox{\r} = (g-1)(\dim(G){-}\rk{G}) + (g'-g)\mbox{\rloc0}
+\sum_{j=1}^2 B_j \mbox{\rloc j}$$
with the local contributions {\rloc j}
recorded in Table~\ref{tab:C}.
(For simplicity of notation, we define $\mbox{\rloc 0} =0$ when there
is no monodromy.)

The exceptional cases are $I_{2k+1}$ with monodromy, and $I_n^*$.  In the
case of $I_{2k+1}$ with monodromy, the formula should be written with a
term $kB_1$ to which the substitution from Proposition \ref{le:su} is {\it
not}\/
applied.\footnote{We are choosing to do this in order to more easily
present the formula as agreeing with a calculation in representation
theory; of course, the version of this formula in which all $B_1$ terms
have been eliminated is also perfectly valid.}

In the case of $I_n^*$, the term $m\bs1\cdot\bso$ in the formula for {\r }
should be broken into two parts, using the substitution from
Proposition \ref{prop:85} to evaluate
a term of the form $(m-2)\bs1\cdot\bso$, but evaluating the remaining term
$2\bs1\cdot\bso$ as
$$2\bs1\cdot\bso=2(-12\ks-m\cc)=(48-8m)(g-1)+(24-2m)B_2$$
(using Proposition \ref{prop:86} for the last step).

The results of all of these manipulations are recorded in the coefficients
given in Table~\ref{tab:C}.

\begin{remark}\label{rem:87}
It is worth observing, for possible generalizations to other cases \cite{AD},
that the substitutions we have used can be formulated intrinsically without
reference to assumptions about the particular types of degenerate fibers
which occur.  This is clear for the substitutions given in 
Propositions~\ref{le:su} and \ref{prop:85}.  In the case of 
Proposition~\ref{prop:86}, when $J$ is finite the substitution only
depends on the discriminant locus.  If $J=\infty$ and we have type
$I_n$ along $\bs1$,  consider the
Weierstrass equation 
\begin{equation}\label{WW} y ^2= x^3 + fx +g \end{equation}
(which is intrinsically associated to the elliptic fibration)
and note that neither $f$ nor $g$ vanishes identically along $\bs1$.
The location of the singularity is given by either $x=-3g/2f$ or
(equivalently) $x=2f^2/9g$.  There is then a divisor $\beta$
on $\bs1$ (in the class $-2\bsi\cdot B$)
represented by $\div_{\bs1}(g)-\div_{\bs1}(f)$ or by
$2\div_{\bs1}(f)-\div_{\bs1}(g)$.  In our case, this divisor coincides
with the divisor $B_1$ (when there is monodromy) or $2B_1$ (when there is
no monodromy) which we used in Proposition~\ref{prop:86}.

Similarly, if $J=\infty$ and we have type $I_n^*$ then neither
$f/s^2$
nor $g/s^3$ vanishes identically along $\bs1$.
The divisor $\beta$ on $\bs1$, which coincides with the divisor
$B_2$ which we used in Proposition~\ref{prop:86}, is represented by
$\div_{\bs1}(g/s^3)-\div_{\bs1}(f/s^2)$ or by
$2\div_{\bs1}(f/s^2)-\div_{\bs1}(g/s^3)$.

Note that this same computation could just as easily be carried out
in the case of multiple components of the discriminant.  The starting
point would be a straightforward generalization of the equation in
Proposition \ref{rev}.  Then for each component of the discriminant,
one would
use the corresponding substitution (according to the singularity type
along that component) and manipulate the substituted formula precisely
as above.  The result is a division into ``non-local'' terms associated
to the various factors of the gauge group (taking precisely the same
form as above), and ``local'' terms associated
to isolated points along the discriminant locus.  We will explore this
generalization further in \cite{AD}.

\end{remark}


\subsection{Step II: Branching rules}\label{sec:bran}

In this subsection and the next, we explain how to systematically determine
representations $\rho_j$, associated to monodromy covers and to degeneration
points, whose charged dimensions reproduce the numbers {\rloc j} which were
calculated in Table~\ref{tab:C}.

Let $\mathfrak h \subset \mathfrak  g$ be a  subalgebra of a Lie algebra.
Given an irreducible
representation $\rho:\mathfrak g\to \GL(N,\mathbb C)$, a natural question
is how $\rho$ decomposes under $\mathfrak h$.
The answer can be obtained by following the ``branching rules''
(see for example \cite{MP}).

\medskip

\noindent \underline{The representation $\rho_0$}

In the case of non-simply laced groups, according to \cite{AKM} the
representation $\rho_0$
is determined by the branching rules for $\mathfrak g_0\subset \mathfrak g$,
where $\mathfrak g_0$ is the non-simply laced algebra and $\mathfrak g$ is
 the corresponding simply laced algebra (whose Dynkin diagram covers that of
 $\mathfrak g_0$).
In each such case, $\mathfrak g_0$ is the fixed subalgebra of some
outer automorphism of $\mathfrak g$ of finite order.

 \begin{proposition}[\cite{MP}] \label{prop96}
The following branching rules hold (using the notation for representations
established in Remark~\ref{rem:notation}):
\begin{itemize}
\item $\sp(k) \subset \su(2k)$ (involutive outer automorphism):
$$\adj \su(2k)=\adj \sp(k) \oplus \Lambda^2_0$$
\item $\sp(k) \subset \su(2k+1)$ (outer automorphism):
$$\adj \su(2k+1)=\adj \sp(k) \oplus \Lambda^2 \oplus \fund
\oplus \fund$$
\item $\so(2k -1) \subset \so(2k)$ (involutive outer automorphism):
$$\adj \so(2k)=\adj \so(2k-1) \oplus \vect $$
\item $G_2\subset \so(8)$ (outer automorphism):
$$\adj \so(8) = \adj G_2 \oplus {\bf 7} \oplus {\bf 7}$$
\item $F_4\subset E_6$ (involutive outer automorphism):
$$ \adj E_6=\adj F_4\oplus{\bf 26}$$
\end{itemize}
\end{proposition}

In the involutive cases, we have $\rho_0$ given as the $(-1)$-eigenspace
of the involution, i.e., the complement of $\adj \mathfrak g_0$ within
$\adj \mathfrak g$.  Thus $\rho_0$
coincides with $\Lambda_0$, $\vect$, and ${\bf 26}$
in the first, third, and fifth cases above, respectively.

In the case of $\sp(k) \subset \su(2k+1)$ , although the automorphism
of $\mathfrak g$ has order $4$, the monodromy action is only order $2$,
and $\rho_0$ is again given by the complement of $\adj \mathfrak g_0$ within
$\adj \mathfrak g$, i.e., $\rho_0=\Lambda^2 \oplus \fund
\oplus \fund$.

In the case of $G_2$, the order $3$ monodromy action leads to the
representation
$\rho_0$ occuring with multiplicity two in the complement of $\adj
\mathfrak g_0$.
(These two copies correspond to the eigenspaces for the monodromy action with
eigenvalues $e^{\pm2\pi i/3}$.)  Thus, in this case $\rho_0={\bf 7}$.

Note that in all cases, the charged dimension of the representation $\rho_0$
agrees with the number {\rloc0} calculated in Table~\ref{tab:C}.

\medskip

\noindent \underline{The representations $\rho_j$}

Representations associated to the points $p$ can also be determined via
branching rules, using a method pioneered by Katz and Vafa \cite{KV}.  If
the general surface section through $p$ has a rational double point
associated to $G'\supset G$, then the representation associated to $p$ is
determined by the corresponding branching rule (modulo a few subtleties to
be discussed in the next subsection).

 \begin{proposition}[\cite{MP}] \label{prop97}
The following branching rules hold (still using the notation from
Remark~\ref{rem:notation}):
\begin{itemize}
\item $\su(n)\subset \su(n+1)$:
$$\adj \su(n+1)=\adj \su(n)\oplus \fund  \oplus \fundbar
\oplus {\bf 1}$$
\item $\su(n) \subset \so (2n)$:
$$\adj \so(2n)=\adj \su (n) \oplus \Lambda^2 \oplus \overline{\Lambda^2}
\oplus {\bf 1}$$
\item $\so(2k)\subset \so(2k+2)$:
$$\adj \so(2k+2) = \adj \so(2k)\oplus \vect  \oplus \vect  \oplus
{\bf 1}$$
\item $\spin(10)\subset E_6$:
$$\adj E_6 = \adj \spin(10) \oplus \spin_+ \oplus \spin_- \oplus {\bf 1}$$
\item $\spin(12)\subset E_7$:
$$\adj E_7 = \adj \spin(12) \oplus \spin_+ \oplus \spin_- \oplus
  {\bf 1} \oplus {\bf 1} \oplus {\bf 1}$$
\item $E_6  \subset  E_7$:
$$\adj E_7= \adj E_6 \oplus {\bf 27}  \oplus \overline{\bf 27}  \oplus {\bf
1}.$$
\item $E_7  \subset  E_8$:
$$\adj E_8= \adj E_7 \oplus {\bf 56}\oplus {\bf 56} \oplus {\bf 1}
\oplus {\bf 1} \oplus {\bf 1}.$$
\end{itemize}
(There are also non-standard embeddings of $D_4$ into $D_5$ which lead
to branching rules involving the $\spin_+$ or $\spin_-$ representations
of $\so(8)$ rather than the vector representation.)
\end{proposition}

Each of these branching rules takes the form
\begin{equation}\label{decomp}
\adj \mathfrak g= \adj \mathfrak g_0 \oplus \rho \oplus \overline{\rho}
\oplus {\bf 1}\end{equation}
for some representation $\rho$; it is $\rho$ which determines the
matter representation.

For example, when the general fiber of type $\su(2k)$, degenerates to
$\su(2k+1)$, then we use the branching rule corresponding
to the inclusion $\su(2k) \subset \su(2k+1)$ to determine
the correct representation ``$\fund$''
 appearing as $\rho$ in the statement of the theorem.

The matter representations $\rho_j$ for non-simply laced groups at non-branch
points can be inferred by looking at the representation of the
corresponding simply laced group.

The cases $\so(12)\subset E_7$ and $E_7  \subset  E_8$ (as well
as the fundamental representation of $\sp(k)$) lead to quaternionic
representations and follow a somewhat different pattern, as we will
explain in the next subsection.  In all other cases, the representation
$\rho_j$ determined by these branching rules has a charged dimension
which agrees with the number {\rloc j} calculated in Table~\ref{tab:C}.

\subsection{Step III: Resolutions of non-generic singularities, deformation
theory,
complex and quaternionic representations}\label{sec:quat}

\quad

Up to this point, we have described the ``matter'' representation as
a complex representation of the group G (as is customary in the
physics literature\footnote{In the physics literature, one refers to
``hypermultiplets taking values in a complex representation'' or,
equivalently, ``half-hypermultiplets taking values in a quaternionic
representation.''}).  However, the representation we need is more
accurately described as a quaternionic representation, that is,
a representation into $\GL(\mathbb{H}^n)$.  Given a complex representation
$\rho$, the representation $\rho\oplus\overline{\rho}$ is automatically
quaternionic---this is how one passes from complex to quaternionic in
many cases.  However, some quaternionic representations cannot be
described as the sum of a complex representation with its complex conjugate.
This explains the presence of the factor $\delta=\frac12$ in certain
terms of the formula for {\r }, since in all cases we are actually
counting $1/2$ of the quaternionic dimension of the representation.

How do these complex and quaternionic representations show up in
the geometry?
Consider again the general elliptic surface passing through $p$.  In all
the cases we are considering, this surface has a rational double point
singularity, which can be associated to a simply laced group $G'$.  Deforming
 to a nearby surface we again find a rational double point, this time the
one associated to the group $G$.

There are three possibilities for a one-parameter family of rational double
points: (1) it fails to admit a simultaneous resolution of singularities,
(2) it is a base-change of a family of type (1) which admits a simultaneous
resolution of singularities, or (3) it admits a simultaneous resolution of
singularities, and is not the base-change of a family which failed to admit
such
a resolution.  When analyzed carefully, the Katz--Vafa prescription \cite{KV}
operates differently in these cases, depending on whether or not
simultaneous resolution is
possible.\footnote{We are grateful to Sheldon
Katz for correspondence on this point.}
It is possible to explicitly compute whether or not this is possible in
each instance, using the
formulas in \cite{KM}.  (One calculates the equation of the family after
performing the base-change which ensures that simultaneous resolution is
possible; the fact that a base-change has been performed can then be
recognized from the dependence of all
coefficients on $t^k$ rather than $t$, for some integer $k$ which
represents the degree of the base-change map.  See \cite{KV}, where many of
these calculations have been carried out.)

Of the branching rules described in Proposition~\ref{prop97}, the first one
($\su(n)\subset \su(n+1)$) falls in case (3),
and all others fall in
cases (1) and (2) (depending on whether $\rho_j$ is being treated
as a representation of a simply laced or a non-simply laced group).  
There is a further distinction that can be made
in case (1): making a base-change to produce a simultaneous resolution,
the base-change group will act on the set of roots, and this action may
or may not induce
monodromy on the Dynkin diagram.

In   case (1), if we perform a finite base-change, a simultaneous
resolution becomes possible and the branching rules determine the
representations which are involved.  However, the covering group for
the base-change acts on these representations, and only the invariant
representation appears in the original family.  In four of the
branching rules from Proposition~\ref{prop97}, there is monodromy
on the Dynkin diagram and we have already analyzed the corresponding
representations from that point of view.  The representation $\rho$
(whose weights are represented by holomorphic curves) is mapped to the
representation $\overline{\rho}$ (whose weights are
represented by anti-holomorphic
curves) with the upshot being that each ramification point on the
parameter curve is associated to $1/2$ of the full representation.
(Of course, we are not counting this as a contribution to the local
representation at the branch point---this part of the representation
theory is non-local, and is accounted for by the representation $\rho_0$.)

Note that these same four branching rules also occur in the context
of case (2) families, where there is no monodromy.  In these cases,
the entire branching rule plays a r\^{o}le, and the quaternionic
representation associated to such a point is $\rho\oplus\overline{\rho}$
(corresponding to the complex representation $\rho$).  Note that
the singularity is fully resolved in these cases, as is reflected in
the Euler characteristic computations in Table~\ref{table:4} in Appendix I.

The remaining two types
of branching rules, $\so(12)\subset E_7$ and $E_7  \subset  E_8$,
only occur in the context of case (1) in our setup, and there is no
monodromy on the Dynkin diagram.
In these cases, the action of the covering group similarly maps
$\rho$ to $\overline{\rho}$, but in these cases the representation
is quaternionic and $\rho\cong\overline{\rho}$.  The upshot is that
the ``complex representation'' associated to each such point is
$1/2$ of the quaternionic representation $\rho$.  (Note that the covering
group acts as $-1$ on the ``${\bf 1}$'' summands in the branching
rule, so that these do not contribute as they are not invariant.)  In both of these cases,
the singularity of the surface is not fully resolved, as is reflected
in the Euler characteristic computations in Table~\ref{table:4} in Appendix I.

The multiplicities of these points are slightly different in the cases
of Kodaira fibers of types II and III, but the same representations
occur.  See \cite{AKM}, where these cases are worked out in detail.

\section{Another look at the substitutions:\\Representation
theory.}\label{sec:sadov}

We have seen how the degenerations
of the general singularity  determine certain representations
of the group $G$; here we show that the converse also holds:
 once one chooses the representations which might occur, the geometry of
the Calabi--Yau
is completely determined by  some relations in representation theory.
We will at the same time verify the additional anomaly cancellations stated
in Theorem~\ref{th:s2}.

We will verify that the generalized Green--Schwarz anomaly cancellation
mechanism works in the way that was proposed by Sadov
\cite{Sadov}.\footnote{We have corrected some minor numerical errors in
\cite{Sadov}.}   The
factored form \eqref{eq:factored} is taken to be
\begin{equation}
\frac12
\left(\frac12 K_B\, \tr R^2 + 2\sum { \bsi} \tr F^2_ i \right)\cdot
\left(\frac12 K_B\, \tr R^2 + 2\sum { \bsi} \tr F^2_ i \right).
\end{equation}
The anomaly cancellation requirements are deduced by comparing this with
equation \eqref{I}.  The coefficients of $(\tr R^2)^2$  agree due to the
relation $9-n_T=K_B^2$, which follows from Noether's theorem on the surface
$B$ (since $\chi(\mathcal{O}_B)=1$).  The remaining coefficients lead to
equations
\begin{align*}
-6\ks (\tr F_i^2)
  &=-\Tr_{\adj} F_i ^2 + \sum_\rho {n_\rho} \Tr _\rho F^2  _i\\
3\cci (\tr F_i^2)^2
  &=-\Tr_{\adj} F_i ^4 + \sum_\rho {n_\rho} \Tr _\rho F^4  _i\\
\bsi\cdot \bsj (\tr F_i^2)(\tr F_j^2)
  &= \sum _{\rho,\sigma} {n_{\rho\sigma}} \ \Tr _{\rho} F^2 _i \  \Tr
_{\sigma} F^2 _j
\end{align*}
which must be evaluated using the relations in the ring of $G$-invariant
functions.  Note that in our case there is a single local factor $G_i$ 
of the gauge group $G$, and we can suppress the subscript $i$
and denote its adjoint curvature by $F$.

We must also specify, for each type of group, a ``fundamental representation''
in which to evaluate the trace $\tr$ on the left-hand side of the equations.
We take $\tr=\Tr_{\fund}$  to be the trace in the
usual fundamental representation for $\su(n)$ and
$\sp(k)$, we take $\tr=\frac12\Tr_{\vect}$ to be one-half of the trace in
the vector representation for $\spin(m)$, and we take $\tr$ to be the
trace in the smallest representation of the group in the case of the
exceptional groups.

Note that if we were to
replace $\tr$ by some multiple of it, say $\lambda \tr$, then we would
multiply $-6K_B\cdot\Sigma_1$ by $\lambda$ and $3\Sigma_1^2$ by
$\lambda^2$.  Making the geometry match the representation theory
completely constrains our choice of $\lambda$, and we express everything
below in terms of the ``correct'' trace for each group.

Having specified the fundamental representation, $\tr F^2$ will correspond to a
basis
of Casimir operators of second order, and $(\tr F^2)^2$ will be one of the
basis elements for
 Casimir operators of the fourth order; when there is a second independent
fourth-order Casimir, the second basis element can be taken to be $\tr F^4$.
Traces taken in other representations can be expressed in terms of these.
We have collected the data of this sort that we need (mostly taken from
Erler \cite{Erler}) in Table~\ref{tab:D} (in which we
use the
notation $\spinrep_*$ to denote either $\spinrep$ or $\spinrep_\pm$).

\begin{table}[ht]
\begin{center}

\begin{tabular}{|c|c|c|c|} \hline
 $G$ & $\rho$ & $\Tr_\rho F^2$ & $\Tr_\rho F^4 $ \\ \hline
 $\su (2) $& $\adj$ &$4\tr F^2$&$8(\tr F^2)^2$ \\
                   & $\fund$ &$\tr F^2$& $\frac12(\tr F^2)^2$ \\ \hline
 $\su (3) $& $\adj$ &$6\tr F^2$&$9(\tr F^2)^2$ \\
                   & $\fund$ &$\tr F^2$& $\frac12(\tr F^2)^2$ \\ \hline
 $\su (n) $,& $\adj$ &$2n\tr F^2$&$6(\tr F^2)^2+2n\tr F^4  $ \\
$n\ge4$         &      $\fund$ &$\tr F^2$& $0(\tr F^2)^2+\tr F^4$ \\
     & $\Lambda^2$ &$(n-2)\tr F^2$&$3(\tr F^2)^2+(n-8)\tr F^4  $  \\
 \hline
 $\sp (k) $,&$\adj$ &$(2k+2)\tr F^2$&$3(\tr F^2)^2+(2k+8)\tr F^4$\\
     $k\ge2$          &$\fund$ &$\tr F^2$&$0(\tr F^2)^2+\tr F^4$\\
     &$\Lambda^2_0$ &$(2k-2)\tr F^2$&$3 (\tr F^2)^2+(2k-8)\tr F^4  $  \\ \hline
$\spin(m)$, & $\adj$ &$(2m-4)\tr F^2$&$12(\tr F^2)^2+(2m-16)\tr F^4  $ \\
$m\ge7$       & $\vect$ &$2\tr F^2$&$0(\tr F^2)^2+2 \tr F^4$ \\
              & $\spinrep_*$ & $\dim(\spinrep_*)(\frac14 \tr F^2)$ &
        $\dim(\spinrep_*)(\frac3{16} (\tr F^2)^2-\frac18 \tr F^4 ) $ \\ \hline
 $E_6$ & $\adj$ &$24\tr F^2$&$18(\tr F^2)^2$ \\
       & $\mathbf{27}$ &$6\tr F^2$&$3 (\tr F^2)^2$\\ \hline
 $E_7$ & $\adj$ &$36\tr F^2$&$24(\tr F^2)^2$\\
       & $\mathbf{56}$ &$12\tr F^2$&$6 (\tr F^2)^2$\\ \hline
 $E_8$ & $\adj$ &$60\tr F^2$&$36(\tr F^2)^2$\\ \hline
 $F_4$ & $\adj$ &$18\tr F^2$&$15(\tr F^2)^2$\\
       & $\mathbf{26}$ &$6\tr F^2$&$3 (\tr F^2)^2$\\ \hline
 $G_2$ & $\adj$ &$8\tr F^2$&$10(\tr F^2)^2$ \\
        & $\mathbf{7}$ &$2\tr F^2$&$(\tr F^2)^2$ \\ \hline
\end{tabular}

\end{center}
\smallskip
\caption{}\label{tab:D}
\end{table}


It is now a straightforward matter to verify the remaining anomaly
cancellations.
We illustrate the procedure in the case of $G=\su(n)$, $n\ge4$, with a matter
representation in which the adjoint representation has multiplicity $g$, 
the fundamental representation has multiplicity $B_2$, and  $\Lambda^2$
has multiplicity $B_1$
  (as specified in Theorem~\ref{th:main}).

{}From Table~\ref{tab:D}, we read off the facts which must hold in order for
the gauge
and mixed anomalies to cancel:
\begin{align*}
-6K_B\cdot\Sigma_1&=2n(g-1)+B_2+(n-2)B_1\\
3\Sigma_1^2&=6(g-1)+0B_2+3B_1\\
0&=2n(g-1)+B_2+(n-8)B_1.
\end{align*}
(Note that there are two equations coming from the quartic anomaly, since
there are two independent fourth order Casimirs.)

To verify these, we use the geometric relations which characterize
$g$, $B_1$, and $B_2$, namely
\begin{align*}
B_1&=-\ks\\
B_2&=(-8K_B-n\bs1)\cdot\bs1\\
g&=(\frac12K_B+\frac12\bs1)\cdot\bs1.
\end{align*}
When these are substituted into the right-hand side of the proposed anomaly
relations,
\begin{align*}
2n(\frac12K_B+\frac12\bs1)+(-8K_B-n\bs1)+(n-2)(-K_B)&=-6K_B\\
6(\frac12K_B+\frac12\bs1)+0(-8K_B-n\bs1)+3(-K_B)&=3\bs1\\
2n(\frac12K_B+\frac12\bs1)+(-8K_B-n\bs1)+(n-8)(-K_B)&=0,
\end{align*}
the relations are verified.

A similar verification can be carried out in all cases.  It is convenient
to supplement the geometric formulas for $g$ and $B_j$'s with a formula for
for $g'-g$ in the case of monodromy, and to compute a quantity $\widehat{B}$
in a few cases (in order to match the representation $\widehat{\rho}$ in the
representation theory, as determined
in Theorem~\ref{th:main}).  We summarize the data in Table~\ref{tab:E}.
(We have omitted the relation $g=(\frac12K_B+\frac12\bs1)\cdot\bs1$, which
always holds.)
Carrying out the verification is then a simple exercise in combining
Tables~\ref{tab:D} and~\ref{tab:E}, as we have done in the case of $\su(n)$
above.

\begin{table}[ht]
{
\renewcommand{\arraystretch}{1}
\begin{center}
\footnotesize
\begin{tabular}{|c|c|c|l|} \hline
Type & $G$ & Basic relations & Derived relations \\ \hline
$I_{2}$ & $\su (2) $& $(-8K_B-2\Sigma_1))\cdot\Sigma_1=B_2$ &\\ \hline
$I_{3}$ & $\su (3) $& $(-9K_B -3\Sigma_1)\cdot\Sigma_1=B_2$ &\\ \hline
$I_{2k} $ & $\sp (k) $, &$-2K_B\cdot\Sigma_1=B_1$ &
         $g'-g=(-\frac12K_B+\frac12\Sigma_1)\cdot\Sigma_1$ \\
             &$k\ge2$& $(-8K_B-2k\Sigma_1)\cdot\Sigma_1=B_2$&\\ \hline
$I_{2k+1} $ & $\sp (k) $, &$-2K_B\cdot\Sigma_1=B_1$ &
         $g'-g=(-\frac12K_B+\frac12\Sigma_1)\cdot\Sigma_1$ \\
             &$k\ge1$& $(-6K_B-(2k+1)\Sigma_1)\cdot\Sigma_1=B_2$&
 $\widehat{B}=2(g'-g)+\frac12B_1+B_2
$\\&&&$\hphantom{\widehat{B}}
=(-8K_B-2k\Sigma_1)\cdot\Sigma_1$\\ \hline
$I_{n}$ & $\su (n) $,& $-2K_B\cdot\Sigma_1=2B_1$ &\\
             &$n\ge4$& $(-8K_B-n\Sigma_1)\cdot\Sigma_1=B_2$ &\\ \hline
$III$ & $\su (2)$ & $(-4K_B-\Sigma_1)\cdot\Sigma_1=B_1$   &
   $\widehat{B}=2B_1=(-8K_B-2\Sigma_1)\cdot\Sigma_1$\\ \hline
$IV$ & $\sp (1) $, &$(-6K_B-2\Sigma_1)\cdot\Sigma_1=B_1$ &
         $g'-g=(-\frac52K_B-\frac12\Sigma_1)\cdot\Sigma_1$ \\
             && &
 $\widehat{B}=2(g'-g)+\frac12B_1
$\\&&&$\hphantom{\widehat{B}}
=(-8K_B-2\Sigma_1)\cdot\Sigma_1$\\ \hline
$IV$ & $\su (3)$ & $(-3K_B-\Sigma_1)\cdot\Sigma_1=B_1$  &
   $\widehat{B}=3B_1=(-9K_B-3\Sigma_1)\cdot\Sigma_1$ \\ \hline
$I_0^*$ & $G_2$ & $(-12K_B-6\Sigma_1)\cdot\Sigma_1=B_1$ &
    $g'-g=(-5K_B-2\Sigma_1)\cdot\Sigma_1$\\ \hline

$I^*_{0}$ & $\spin (7) $& $(-4K_B-2\Sigma_1)\cdot\Sigma_1=B_1$   &
       $g'-g=(-\frac32K_B-\frac12\Sigma_1)\cdot\Sigma_1$ \\
&&$(-4K_B-2\Sigma_1)\cdot\Sigma_1=B_2$&\\ \hline
$I^*_{0} $ & $\spin (8) $&$(-2K_B-\Sigma_1)\cdot\Sigma_1=B_1$  &\\
   &&$(-4K_B-2\Sigma_1)\cdot\Sigma_1=B_2$&    \\ \hline

$I^*_{1}$ & $\spin (9) $& $(-6K_B-4\Sigma_1)\cdot\Sigma_1=B_1$   &
        $g'-g=(-\frac52K_B-\frac32\Sigma_1)\cdot\Sigma_1$ \\
&&$(-2K_B-\Sigma_1)\cdot\Sigma_1=B_2$&\\ \hline
$I^*_{1} $ & $\spin (10) $&$(-3K_B-2\Sigma_1)\cdot\Sigma_1=B_1$  &
        \\  &&$(-2K_B-\Sigma_1)\cdot\Sigma_1=B_2$&\\ \hline
$I^*_{2}$ & $\spin (11) $& $(-8K_B-6\Sigma_1)\cdot\Sigma_1=B_1$   &
       $g'-g=(-\frac72K_B-\frac52\Sigma_1)\cdot\Sigma_1$ \\
&&$(-2K_B-\Sigma_1)\cdot\Sigma_1=B_2$&\\ \hline
$I^*_{2} $ & $\spin (12) $&$(-4K_B-3\Sigma_1)\cdot\Sigma_1=B_1$  &\\
   &&$(-2K_B-\Sigma_1)\cdot\Sigma_1=B_2$&    \\ \hline

$I^*_{2k-3}$ & $\so (4k +1) $& $(-6K_B-2k\Sigma_1)\cdot\Sigma_1=B_1$   &
        $g'-g=(-\frac52K_B-(k{-}\frac12)\Sigma_1)\cdot\Sigma_1$ \\
&&$(-2K_B-\Sigma_1)\cdot\Sigma_1=0$&\\ \hline
$I^*_{2k-3} $ & $\so (4k+2) $&$(-3K_B-k\Sigma_1)\cdot\Sigma_1=B_1$  &
        \\  &&$(-2K_B-\Sigma_1)\cdot\Sigma_1=0$&\\ \hline
$I^*_{2k-2}$ & $\so (4k +3) $& $(-8K_B-(2k+2)\Sigma_1)\cdot\Sigma_1=B_1$   &
       $g'-g=(-\frac72K_B-(k{+}\frac12)\Sigma_1)\cdot\Sigma_1$ \\
&&$(-2K_B-\Sigma_1)\cdot\Sigma_1=0$&\\ \hline
$I^*_{2k-2} $ & $\so (4k+4) $&$(-4K_B-(k+1)\Sigma_1)\cdot\Sigma_1=B_1$  &\\
   &&$(-2K_B-\Sigma_1)\cdot\Sigma_1=0$&    \\ \hline
 $IV ^*$ & $E_6$ & $(-12K_B-8\Sigma_1)\cdot\Sigma_1=4B_1$ &\\ \hline
$IV ^*$ & $F_4$ &  $(-12K_B-8\Sigma_1)\cdot\Sigma_1=2B_1$ &
     $g'-g=(-\frac52K_B-\frac32\Sigma_1)\cdot\Sigma_1$    \\ \hline
 $III ^*$ & $E_7$ & $(-12K_B-9\Sigma_1)\cdot\Sigma_1=3B_1$ &\\ \hline
 $II ^*$ & $E_8$ & $(-12K_B-10\Sigma_1)\cdot\Sigma_1=0$   &\\ \hline
\end{tabular}
\end{center}
}
\smallskip
\caption{The relations.}\label{tab:E}

\end{table}


\begin{remark}\label{rem:final}
As in Remark~\ref{rem:87}, we can express this part of the verification
of the anomaly cancellation in terms which are somewhat more intrinsic.
We defer the details of this to \cite{AD}, but observe here how this
can be carried out in the case of $\su(n)$, $n\ge4$.  

The intrinsic geometric quantities we need are the divisor $\bso\cdot\bs1$,
the arithmetic genus $p_a(\bs1)$, and the divisor $\beta$ from 
Remark~\ref{rem:87} (which is the intrinsic version of $2B_1$).
We derive from these an intrinsic version of $B_2$, represented as 
$\bso\cdot\bs1-2\beta$.  Then in the anomaly cancellation requirements,
we can represent the coefficient of $\tr F_i^2$ as
$$ 2n(p_a(\bs1)-1)+(\bso\cdot\bs1-2\beta)+\frac{n-2}2\beta$$
and the coefficient of $(\tr F_i^2)^2$ as
$$ 6(p_a(\bs1)-1)+\frac{3}2\beta.$$

\end{remark}

\newpage

\section*{Appendix I: How to compute {\r}\\
 (the coefficients in Proposition \ref{rev} and other things)}
\label{sec:tables}\renewcommand{\thesection}{I}\setcounter{theorem}{0}
\setcounter{subsection}{0}\setcounter{table}{0}
\makeatletter
\renewcommand\thetable{\@arabic\c@table}
\makeatother

In this section we study the local equations and the  geometric data
for each group and their generic degenerations.

Following  \cite{BIKMSV} we analyze the local equations in 
Tables~\ref{table:1} and \ref{table:2}.
  In Tables~\ref{table:3} and \ref{table:4}
we list, for each group, the coefficients of the right
hand side of the equation defining {\r}, in Proposition \ref{rev}.
The entries of Table~\ref{table:1} are taken from
\cite{BIKMSV}, those of Table~\ref{table:3} are well known; to compute the
others
we need the affine equations of \eqref{Wg} and \eqref{fg}.
We will work out the details for the case $G=\su(2k)$ in Appendix II.

We  need to use a more
general form of the Weierstrass equation \eqref{W}, namely
\begin{equation} \label{Wg} y^2 + a_1xy+ a_3 y= x^3 + a_2x^2 + a_4x + a_6
.\end{equation}

 Since $W$ is assumed to be a Calabi--Yau $a_j \in |-j\kb|$.

\begin{definition}\label{fg} It is convenient to use the following:

$b_2 = {a_1} ^2 + 4a_2$, \ $b_4=a_1 a_3 + 2a_4$, \ $b_6= a_3 ^2 +4a_6 $,

$b_8= ({a_1}^2 + a_2 ) a_6  -{a_4}^2 $

The coefficients in \eqref{W} are now:
\begin{equation*}f=\frac{-1}{48}({b_2} ^2 - 24 b_4), \  \
g=\frac{-1}{864}({-b_2} ^3 + 36
b_2
b_4 -216 b_6).\end{equation*}

If $a_j$ (resp.\ $b_j$) vanish along $\bs1$ of order $k$,
then we write
$$a_{j,k}=\frac{a_j}{s^k} \ (\text{resp. } \ b_{j,k}=\frac{b_j}{s^k}).$$
\end{definition}

Table~\ref{table:1} is mostly taken from \cite{BIKMSV}:
the first two columns list the Kodaira fiber and the associated group
(see Section \ref{sec:G}); in the middle columns we write the order
of vanishing of each $a_i$ along $\bs1$.
Recall that our hypothesis (a flat Calabi--Yau fibration) imposes
some restriction on the self-intersection of the ramification
divisor (see the Remark after the Main Theorem \ref{th:main}).
In the last column, we exhibit how the equation for $\bso \bmod s$ breaks
into factors; the power $r_j$ which gives the multiplicity of the factor
$\beta_j$ is indicated in the factorization in each case.

We have incorporated some necessary corrections to the Table from
\cite{BIKMSV}.
First,  the
entry for
$I_{2k+1}$, $k\ge1$, with gauge group $\su(2k{+}1)$ corresponds to the
Weierstrass equation
$$y^2+a_1xy+a_{3,k}s^ky=
x^3+a_{2,1}sx^2+a_{4,k+1}s^{k+1}x+a_{6,2k+1}s^{2k+1},$$
which has discriminant
\[
 - {\displaystyle \frac {1}{16}} \,{a_{1}}^{4}\,({a_{1}}^{2}\,{a
_{6, \,2\,k + 1}} - {a_{1}}\,{a_{3, \,k}}\,{a_{4, \,k + 1}} + {a
_{3, \,k}}^{2}\,{a_{2, \,1}})
s^{2k+1}  -\frac1{16}a_1^3a_{3,k}^3s^{3k}+ O(s^{2k+2}).
\]
Thus, the correct  leading
term in the local equation of $\bso$ in this case (the ``residual
discriminant'') takes the form
\[a_1^3(a_1b_{8,2k+1}-a_{3,k}^3), \quad \text{if}\ k=1,\]
and
\[a_1^4b_{8,2k+1},\quad \text{if}\ k>1\]
({\it not}\/ $a_1^6 a_{6,2k+1}$ as was written in
\cite{BIKMSV}).

Second, the residual discriminant in the case $IV$ (with gauge group
$\su (3)$) should read
$-27a_{3,1}^4$ rather than $-27a_{3,2}^4$.

\begin{table}[ht]

\begin{center}
\begin{tabular}{|c|c|c|c|c|c|c|c|} \hline
Type & $G$ & $a_1$ & $a_2$ & $a_3$ & $a_4$ & $a_6$ & $(\beta_1)^{r_1}
  (\beta_2)^{r_2}$ \\ \hline
$I_{1}$ & $\{e\}$& $0$ & $0$ & $1$ & $1$ & $1$ & $(b_2)^3 (a_{6,1})^1$ \\
  \hline
$I_{2}$ & $\su(2)$ & $0$ & $0$ & $1$ & $1$ & $2$ & $(b_2)^2 (b_{8,2})^1$ \\
  \hline
$I_{3}$ & $\su(3)$ & $0$ & $1$ & $1$ & $2$ & $3$ & $(a_1)^3
  (a_1b_{8,3}-a_{3,1}^3)^1$ \\ \hline
$I_{2k}, k \geq 2$ & $\sp(k)$ & $0$ & $0$ & $k$ & $k$ & $2k$ & $(b_2)^2
  (b_{8,2k})^1$ \\ \hline
$I_{2k+1}, k \geq 1$ & $\sp(k)$ & $0$ & $0$ & $k{+}1$ & $k{+}1$ & $2k{+}1$
  & $(b_2)^3 (a_{6,2k+1})^1$ \\ \hline
$I_n, n\geq4$ & $\su(n)$ & $0$ & $1$ & $[\frac n2]$ & $[\frac{n{+}1}2]$ &
  $n$ & $(a_1)^4 (b_{8,n})^1$ \\ \hline
$II$ & $\{e\}$ & $1$ & $1$ & $1$ & $1$ & $1$ & $(a_{6,1})^2$  \\ \hline
$III$ & $\su(2)$ & $1$ & $1$ & $1$ & $1$ & $2$ & $(a_{4,1})^3$  \\ \hline
$IV$ & $\sp(1)$ & $1$ & $1$ & $1$ & $2$ & $2$ & $(b_{6,2})^2$  \\ \hline
$IV$ & $\su (3)$ & $1$ & $1$ & $1$ & $2$ & $3$ & $(a_{3,1})^4$  \\ \hline
$I_0^*$ & $G_2$ & $1$ & $1$ & $2$ & $2$ & $3$ & $(\Delta_{12,6})^1$  \\ \hline
$I_0^*$ & $\spin(7)$ & $1$ & $1$ & $2$ & $2$ & $4$ & $(a_{2,1}^2-a_{4,2})^1
  (a_{4,2})^2$ \\ \hline
$I_0^*$ & $\spin(8)$ & $1$ & $1$ & $2$ & $2$ & $4$ &
  $\left(\sqrt{a_{2,1}^2-a_{4,2}}\right)^2 (a_{4,2})^2$ \\ \hline
$I_1^*$ & $\spin (9)$ & $1$ & $1$ & $2$ & $3$ & $4$ & $(b_{6,4})^1 (a_{2,1})^3$
  \\ \hline
$I_1^*$ & $\spin(10)$ & $1$ & $1$ & $2$ & $3$ & $5$ & $(a_{3,2})^2 (a_{2,1})^3$
  \\ \hline
$I_2^*$ & $\spin(11)$ & $1$ & $1$ & $3$ & $3$ & $5$ &
  $(a_{4,3}^2-4a_{2,1}a_{6,5})^1 (a_{2,1})^2$  \\ \hline
$I_2^*$ & $\spin (12)$ &1& $1$ & $3$ & $3$ & $5$ &
  $\left(\sqrt{a_{4,3}^2-4a_{2,1}a_{6,5}}\right)^2 (a_{2,1})^2$ \\ \hline
$I^*_{2k-3}, k\geq3 $ & $\so(4k{+}1)$ & $1$ & $1$ & $k$ & $k{+}1$ & $2k$ &
  $(b_{6,2k})^1 (a_{2,1})^3$ \\ \hline
$I^*_{2k-3}, k\geq3$ & $\so({4k{+}2})$ & $1$ & $1$ & $k$ & $k{+}1$ &
  $2k{+}1$ & $(a_{3,k})^2 (a_{2,1})^3$ \\ \hline
$I^*_{2k-2}, k\geq3 $ & $\so (4k{+}3)$ & $1$ & $1$ & $k{+}1$ & $k{+}1$ &
  $2k{+}1$ & $(a_{4,k+1}^2{-}4a_{2,1}a_{6,2k+1})^1 (a_{2,1})^2$ \\ \hline
$I^*_{2k-2}, k\geq3$ & $\so(4k{+}4)$ & $1$ & $1$ & $k{+}1$ & $k{+}1$ &
  $2k{+}1$ & $\left(\sqrt{a_{4,k+1}^2{-}4a_{2,1}a_{6,2k{+}1}}\right)^2
  (a_{2,1})^2$ \\   \hline
$IV ^*$ & $F_4$ & $1$ & $2$ & $2$ & $3$ & $4$ & $(b_{6,4})^2$  \\ \hline
$IV ^*$ & $E_6$ & $1$ & $2$ & $2$ & $3$ & $5$ & $(a_{3,2})^4$  \\ \hline
$III ^*$ & $E_7$ & $1$ & $2$ & $3$ & $3$ & $5$ & $(a_{4,3})^3$  \\ \hline
$II ^*$ & $E_8$ & $1$ & $2$ & $3$ & $4$ & $5$ & $(a_{6,5})^2$  \\ \hline
\end{tabular}
\end{center}
\smallskip
\caption{}\label{table:1}
\end{table}

\newpage

\begin{remark}
Following \cite{BIKMSV} (see also Proposition \ref{link}) we see that
 if
$$\bso \cap \bs1= \{ P_1 ^1, \cdots P^{B_1} _1, P^1 _2, \cdots P^{B_2} _2
\},$$
then the local equation of $\bso$  around $P^\ell _i$
does not depend on $\ell$, but only on $i=1,2$.
\end{remark}

In  Table~\ref{table:2} we list the
 local equation (l.e.) of $\bso$ around $P_1$ and $P_2$.
As usual, we denote by $s=0$
the divisor $\bs1$; $t$ is a convenient coordinate vanishing at
$P_i$ and $\g _i$ is a suitable invertible function near $\{s=t=0\}$.

Our assumption on the existence of a smooth Calabi--Yau resolution
 imposes
of $\bs1$ and $\bso$
 We write ``\nonmin" or ``\terminal''
if the intersection type, as stated in Table~\ref{table:1} is
not compatible with our hypothesis
due to the singularities being non-minimal or having no small resolution.

\begin{table}[ht]

\begin{center}

\begin{tabular}{|c|c|c|c|} \hline
 Type & $ G $  & l.e. at $P_1$  & l.e. at $P_2$ \\ \hline
$I_{1}$ & $\{e\}$& $\gamma _1 t^3 + \gamma _2 s=0$ & \terminal \\ \hline
$I_{2}$ & $\su (2) $& $\gamma _0 t^2 + \gamma _1 s =0$&transversal \\ \hline
$I_{3}$ & $\su (3) $& $\g _0 t^3 + \g _1 s=0$&  transversal
  \\ \hline
$I_{2k}, k \geq 2$ & $\sp (k) $ & $t^2 - \gamma s^k=0$& transversal  \\
\hline
$I_{2k+1}, k \geq 2$ & $\sp (k) $&$t^2(\g _0 t + \g
_1s) + \g _2 t s^{k+1} +
\g_3 s ^{k+2}=0$&    \terminal  \\ \hline
$I_{n}, n \geq 4$ & $\su (n) $   & $s^{n} + \gamma _1 t^2 = 0$&
transversal \\ \hline
$II$ & $\{e\}$ & \terminal &  \\ \hline
$III$ & $\su (2)$ & $\gamma _0 s + \gamma _1 t^3=0$ &  \\ \hline
$IV$ & $\sp (1)  $ & $\g t^2 + s=0$ & \\ \hline
$IV$ & $\su (3) $ & $t^4 + \gamma _0 s^2 + \gamma _1 s t^2=0$ & \\ \hline
$I_0^*$ & $G_2$ &   transversal &\\ \hline
$I_0^*$ & $\spin (7)$ & transversal & $\g t^2 +s=0$  \\ \hline
$I_0^*$ & $\spin (8)$ &$\g _0 s + \g _1 t^2=0$& $\g _0 s + \g _1 t^2$ \\ \hline
$I^*_{1}$ & $\spin (9)$ & transversal & $\g _0 s + \g _1 t^3$\\ \hline
$I^*_{1}$& $\spin ({10})$  & $\g t^2 +s=0$&$\g t^2 +s=0$\\ \hline
$I^*_{2}$ & $\spin (11)$ & transversal & $\g t^2 +s=0$ \\ \hline
$I^*_{2}$ & $\spin(12)$ &  $\g t^2 +s=0$ & $\g t^2 +s=0$\\ \hline
$I^*_{n}, n\geq3$ & $\so (2n+7)$ & transversal & \nonmin\\ \hline
$I^*_{n}, n\geq3$& $\so (2n+8)$  & $\g t^2 +s=0$&\nonmin\\ \hline
 $IV ^*$ & $F_4$ & $\g s + t^2=0$ &  \\ \hline
 $IV ^*$ & $E_6$ & $\gamma _0 s + t^4=0$ & \\ \hline
$III ^*$ & $E_7$ & $\gamma _0 s + \gamma _ 1 t^3=0$ &  \\ \hline
$II^*$ & $E_8$ & \nonmin & \\ \hline
\end{tabular}
\end{center}
\smallskip
\caption{}\label{table:2}
\end{table}


In Table~\ref{table:3}, $h$ denotes the Coxeter number of the group $G$,
$m$ the multiplicity of $\bs1$ in the discriminant,
and  $\mu(f)$ (resp.\ $\mu(g)$) the vanishing of $f$ (resp.\ $g$)
in equation \eqref{W}
along $\bs1$ (see also Section \ref{sec:wa}).

\begin{table}[ht]
{\renewcommand{\arraystretch}{1.5}

\begin{center}

\begin{tabular}{|c|c|c|c|c|c|c|} \hline
 Type & $ G $& h &{\boldmath  $\rk$} & $m$ & $ \mu (f) $ &  $\mu (g)$ \\
\hline
$I_{1}$ & $\{e\} $&  -- & -- &   $1$ &$ 0 $ &  $0$ \\ \hline
$I_{2}$ & $\su (2) $&  $2$ & $1$ &  $2$  & $0$ & $0$  \\ \hline
$I_{3}$ & $\su (3) $& $3$ &$2$ & $3$ & $0$ &$0$\\ \hline
$I_{2k}, k \geq 2$ & $\sp (k) $&  $2k$ & $k$ &   $2k$ & $ 0 $ &  $0$\\ \hline
$I_{2k+1}, k \geq 1$ & $\sp (k) $&  $2k$ & $k$ &   $2k+1$ & $ 0 $ &  $0$\\
\hline
$I_{n}, n \geq 4$ & $\su (n) $&  $n$ & $n-1$ &   $n$ & $ 0 $ &  $0$\\
\hline
$II$ & $\{e\}$& -- & -- & $2$ & $1$ & $1$ \\ \hline
$III$ & $\su (2)$ & $2$ & $1$ & $3$ & $1$ & $2$  \\ \hline
$IV$ & $\sp (1)  $ & $2$ & $1$& $4$ & $2$ & $2$\\ \hline
$IV$ & $\su (3)$ & $3$ & $2$ & $4$ & $2$ & $2$  \\ \hline
$I_0^*$ & $G_2$ &  $6$& $2$& $6$ & $2$ &$3$ \\ \hline
$I^*_{n}, n\ge0$ & $SO(2n+7)$ & $2n+6$& $n+3$& $n+6$ &$2$  & $3$\\ \hline
$I^*_{n}, n\ge0$ & $SO(2n+8)$ & $2n+6$& $n+4$& $n+6$ &$2$  & $3$\\ \hline
 $IV ^*$ & $F_4$ & $12$ & $4$& $8$& $3$ &  $4$  \\ \hline
 $IV ^*$ & $E_6$ & $12$ & $6$ & $8$  & $3$ & $4$ \\ \hline
 $III ^*$ & $E_7$ & $18$ & $7$ & $9$ & $3$ & $5$ \\ \hline
 $II ^*$ & $E_8$ & $30$ & $8$ & $10$ & $4$ & $5$ \\ \hline

\end{tabular}

\end{center}
}
\smallskip
\caption{}\label{table:3}
\end{table}


 In  Table~\ref{table:4} we write, for each Kodaira type fiber and associated
group, the coefficients needed to compute {\r}, as in Proposition
\ref{rev}.
The general Kodaira type fiber over $\bs1$ degenerates over
both $P_i$ at the intersection with $\bso$.
As in Table~\ref{table:2} we write ``\nonmin" or ``\terminal'' if the intersection
type, as stated in Table~\ref{table:1} is
not compatible with our hypothesis.
We describe the degenerate singular fibers:
if they are of Kodaira type we use Kodaira's notation.
Note that these are not necessarily the Kodaira type of the general
Weierstrass surface passing through the degenerate fiber;
for example in the case of $G=E_7$ ($III^*$), the degenerate fiber is
again of type $III^*$, but the general Weierstrass surface has
a $II^*$ singularity (see also Section \ref{sec:quat}).
These distinctions are important in computing {\r} as in Theorem
\ref{th:main}.


The fibers of non-Kodaira type are the branch points of an
outer automorphism of
the group; we denote these with ``br.".

\begin{table}[ht]
{
\renewcommand{\arraystretch}{1.5}

\begin{center}
\begin{tabular}{|c|c|c|c|c|c|c|c|} \hline
 Type & $ G $  & $\mu_1 (f,g) $ & $\mu_2(f,g)$ &
$ \epsilon _1$ & $\epsilon _2$  & $\c (X_{P_1})$ & $ \c(X_{P_2}) $\\ \hline
$I_{1}$ & $\{e\} $& $2$ & $0$ &
$ -1$ & $-1$   & $2$ ($II$) & \terminal \\ \hline
$I_{2}$ & $\su (2) $& $3$ &$0$ &
  $ -1$ & $-1$  & $3$ ($III$) &  $3$ ($I_3$)  \\ \hline
$I_{3}$ & $\su (3) $& $8$ &$0$ &
$ -1$ & $-1$ &  $4$ ($IV$)& $ 4 $ ($I_4$) \\ \hline
$I_{2k}, k \geq 2$ & $\sp (k) $&  $3k$ & $0$ &  $ k-2 $ & $-1$
&    $k+2$ (br.)& $2k+1$ ($I_{2k+1}$)  \\ \hline
$I_{2k+1}, k \ge1$ & $\sp (k)$ & $3k+3$ & $0$&  $k+2$& $-1$  & $k+2$
(br.)& \terminal    \\ \hline
$I_{n}, n \geq 4$ & $\su (n) $& $3n$ & $0$ &  $n-2$ & $-1$
  & $n+2 $ ($D_{n}$)& $n +1$  ($I_{n+1}$) \\ \hline
$II$ & $\{e\}$ &0 & &$-1$  &  &  \terminal &\\ \hline
$III$ & $\su (2)$ & $0$ & &    $-1$ &  &  $4$ $(IV)$  &\\ \hline
 $IV$ & $\sp (1)  $ & $0$ & &   $-1$ &  & $3$        (br.) &\\ \hline
$IV$ & $\su (3) $ & $0$ & &    $2$ & & $6$  ($I_0^*$) &\\ \hline
$I_0^*$ & $G_2$ & $0$ & &    $-1$& & $5$  (br.) & \\ \hline
$I_0^*$ & $\spin (7)$ & $0$ &$0$ &    $-1$& $-1$&  $5$ (br.)& 7
($I^*_1$)
 \\ \hline
$I_0^*$ & $\spin (8)$ &$0$&$0$ &  $-1$& $-1$& $7 \ (I^*_1)$& $7$ ($I^*_1$) \\
   \hline
$I_1^*$ & $\spin (9)$ &$0$& $2$&   $-1$ &$-1$ &6 (br.) &  $8$ ($IV^*$)\\
\hline
$I_1^*$ & $\spin (10)$ &$0$& $2$&  $-1$ & $-1$& $8 \ (I^*_2)$ &
  $8$ ($IV^*$) \\ \hline
$I_2^*$ & $\spin (11)$ &$0$& $3$&   $-1$&$-1$ & $7$ (br.)& $8$ ($I_2^*$) \\
\hline
$I_2^*$ & $\spin (12)$ &$0$& $3$&   $-1$&$-1$ &$9 \ (I^*_3)$ &
  $8$ ($I_2^*$) \\ \hline
$I^*_{n}, n\ge3 $& $\so(2n+7)$ &$0$& \nonmin   &
  $-1$&\nonmin& $n+5$ (br.)  &\nonmin\\ \hline
$I^*_{n}, n\ge3$& $\so({2n+8})$ &$0$& \nonmin &     $-1$
   &\nonmin&  $n+7$ ($I^*_{n+1}$) &\nonmin\\ \hline
 $IV ^*$ & $F_4$ & $0$ & &     $-1$ &  &  $6$ (br.)&\\ \hline
 $IV ^*$ & $E_6$ & $0$ & &     $-1$  &  &$9$ ($III^*$)&\\ \hline
$III ^*$ & $E_7$ & $0$ & &   $-1$ &    & $9$ ($III^*$) &\\ \hline
$II ^*$ & $E_8$ & \nonmin & &\nonmin &    &\nonmin  &\\ \hline
\end{tabular}
\end{center}
}
\smallskip
\caption{}\label{table:4}
\end{table}


\section*{Appendix II: The entries in the above Tables for $G= \su(2k), \ k
\geq 2$ and $I_{2k}$ fiber type.}
\label{sec:entries}\renewcommand{\thesection}{II}\setcounter{theorem}{0}
\setcounter{subsection}{0}

We illustrate the pattern of computations needed to compile the Tables in
Appendix I with the specific example $G=\su(2k)$.

The generalized Weierstrass equation has the form:
$$y^2 + a_1xy= x^3 + a_2sx^2 + a_4s^k x+ a_6 s^{2k}.$$

\smallskip

$b_2 = {a_1} ^2 + 4a_2s$, \ $b_4= 2a_4s^k$, \ $b_6= 4a_6 s^{2k}$,

\smallskip

$b_8= [({a_1}^2 + 4a_2 s) a_6  -{a_4}^2] s ^{2k}, \quad b_{8,2k}= ({a_1}^2
+ 4 a_2 s) a_6  -{a_4}^2$

\smallskip

$f=\frac{-1}{48}({b_2} ^2 - 24 b_4), \quad  \ g=\frac{-1}{864}({-b_2} ^3 + 36
b_2
b_4 -216 b_6).$

\smallskip

$$\aligned \bsa : \phantom{xx}  &s^{2k}\{-({a_1 }^4+ 16 {a_2 }^2 s^2 + 8
{a_1
}^2 a_2 s)[({a_1 }^2 + 4 a_2 s) a_6 - {a_4 }^2] +\\
& - 8s ^k [8 {a_4 }^3 +  27 \cdot 2{a_6
}^2s^k- 9 a_4a_6 ({a_1 }^2 + 4 a_2s)],
\endaligned
$$

$$\aligned  \bso : \phantom{xx} &-({a_1 }^4+ 16 {a_2 }^2 s^2 + 8 {a_1 }^2
a_2
s)[({a_1 }^2 + 4 a_2 s) a_6 - {a_4 }^2] +\\
& - 8s ^k [8 {a_4 }^3 +  27 \cdot 2{a_6
}^2s^k- 9 a_4a_6 ({a_1 }^2 + 4 a_2s)]\}.
\endaligned
$$

At the points of intersections of $\bso$ and $\bs1$, either $a_1=0$
($P^{\ell} _1$)
or $b_{8,2k}=0$ ($P^\ell _2$).  (In the notation of Section
\ref{sec:wa}, $a_1=\beta_1$.)
\begin{remark}$r_1=4$ and $r_2=1$; there are $B_1= -\ks$
points of $P_1$ type,
and $(-8\kb -2k \bs1)\cdot \bs1=B_2$ points of $P_2$ type.
The second condition follows from the first one, as
$\bs1 \cdot \bso = 4B_1 + B_1$.
\end{remark}

\subsection{Computing $\epsilon _1$:} \label{subsec:2b}

Let $t=: a_1, s$ be the local coordinates around a point $P^1_ \ell$. (In
the notation of Section
\ref{sec:wa}, $a_1=\beta_1$.)

Then

$$\bso : \gamma _0 t^4 + \gamma _1 s^2 + \gamma _2 t^2 s + \gamma _3 s^k,$$
where $\gamma _i$ is invertible at $s=t=0$.

We can write
$$\bso : \gamma _0 t^{2k} +  \gamma _1 s^2 =0$$
which defines an $A_{2k-1}$ curve singularity. Since the blowup of an
$A_{2k-1}$ curve singularity
yields an $A_{2k-3}$ singularity, we have

$$(\# \phi ^{-1}(P_1);\ \{ \alpha ^1 _v \} )= (2;2, \cdots 2), \ (k \text{
times});\ \ \text{ then }\epsilon _1 = 2k-2.$$

\subsection{Computing $\epsilon _2 $:} \label{subsec:2a}

Since $\bso$ is smooth around each point $P_2$

$$(\# \phi ^{-1}(P_2);\ \{ \alpha ^2 _v \} )= (1;1); \epsilon _2 = -1.$$

\subsection{Computing $\mu (f,g)$: } \label{subsec:2c}

{}From the equations we see that $f$ and $g$ have a common zero along $\bs1$
when $b_2=0$, and there $-2 \kb \cdot \bs1$ such points.
Now set
 $$g' = \frac{b_2}{18} f + g= \frac{1}{72}(-b_2 b_4)
-  \frac{1}{12} b_6 .$$
Then $\mu (f,g) = \mu (f, g')$ \cite[Section 1]{F}.

{}From the equation  above
we see that $P \in \bs1$ is a common zero of $f$ and $g$ if and only if
$a_1=0$.
As in \ref{subsec:2b} we take $t:= a_1, s$ as the local coordinates around
$P$.

$$\aligned
\mu (f, g) &= \dim _{\mathbb C} \mathbb C [[s,t]] / (f, g'), \text{ where}\\
 f &\approx t^4 + \gamma _1 t^2 s + \gamma _2 s^2  + \gamma _3  s^k \\
 g' &\approx 24 s^k \{ \gamma _4t +  s^{k})\},
\endaligned$$
for suitable invertible functions $\gamma _i$ (around $s=t=0$).
 Then  \cite[Ex. 1.2.5]{F}
$$\mu (f, g)=  6k . $$

\subsection{Computing  $\c (X_{P_1})$.}\label{subsec:degf}
After $\ell$ blowups the Weierstrass equation becomes:
$y ^2 + a_1 xy = x^3 s^{\ell} + a_2 s x^2 + a_4 x s^{k - \ell}+ a_6s ^{2k
-2\ell}$, and there are isolated singular points (nodes) on the
fiber at $P_1=0$. These points can be blown up with small resolutions:
the fiber over the points $P_1$ is of Kodaira type $D_{2k}$ and $\c
(X_{P_1})= 2k+2.$

\subsection{Computing  $\c (X_{P_2})$.}\label{subsec:degff}
$\bso$ and $\bs1$ intersect transversally at $P_2$, and it is easy to see
that the corresponding fiber
$X_{P_2}$ is of type $I_{2k+1}$ and $\c (X_{P_2})= 2k+1.$


\begin{thebibliography}{10}

\bibitem{Asp}
P.~S. Aspinwall and M.~Gross, {\em The {$SO(32)$} heterotic string on a {K3}
  surface}, Phys. Lett. B {\bf 387} (1996) 735--742, {\urlfont
hep-th/9605131}.

\bibitem{AKM}
P.~S.~Aspinwall, S.~Katz, and D.~R.~Morrison, {\em Lie groups, Calabi--Yau
threefolds and $F$-theory}, {\urlfont hep-th/0002012v2}.

\bibitem{BKKM} P.~Berglund, S.~Katz, A. Klemm, and P. Mayr,
{\em New Higgs transitions between dual $N{=}2$ string models},
 Nucl. Phys. B {\bf 483} (1997) 209--228, {\urlfont hep-th/9605154}.

\bibitem{BIKMSV}
M.~Bershadsky, K.~Intriligator, S.~Kachru, D.~R.~Morrison, V.~Sadov, and
C.~Vafa, {\em
  Geometric singularities and enhanced gauge symmetries}, Nucl. Phys. B {\bf
  481} (1996) 215--252, {\urlfont hep-th/9605200}.

\bibitem{Bri:Nice}
E. Brieskorn, {\em Singular elements of semi-simple algebraic groups},
{A}ctes
  {C}ongr{\`e}s intern. {M}ath. 1970, tome 2, Gauthier-Villars, Paris, 1971,
  pp.~279--284.

\bibitem{CPR}
P.~Candelas, E.~Perevalov, and G.~Rajesh,
\newblock {\em Matter from toric geometry},
\newblock Nucl. Phys. {\bf B519} (1998) 225--238, {\urlfont hep-th/9707049}.

\bibitem{coxeter:annals}
H. S. M. Coxeter, {\em Discrete groups generated by reflections}, Annals of
Math. (2) {\bf 35} (1934) 588--621.

\bibitem{coxeter:weyl}
H. S. M. Coxeter, {\em Discrete groups generated by reflections}, in: ``The
Structure and Representation of Continuous Groups,'' by Hermann Weyl,
Institute for Advanced Study, 1935, pp.186--210.

\bibitem{deligne} P.~Deligne, {\em Le s\'erie exceptionnelle de groupes de
Lie},
C.R. Acad. Sci. Paris {\bf t. 322, S\'erie I} (1996) 321--326.

\bibitem{DE}
D.-E. Diaconescu and R.~Entin,
\newblock {\em Calabi--Yau spaces and five-dimensional field theories with
  exceptional gauge symmetry},
\newblock Nucl. Phys. {\bf B538} (1999) 451--484, {\urlfont hep-th/9807170}.

\bibitem{DV} P.~Du Val, {\em On isolated singularities which do not affect
the condition of adjunction, Part I},
Proc. Cambridge Phil. Soc {\bf  30} (1934) 453--465.

\bibitem{DV:book} P. Du Val, ``Homographies, Quaternions, and
Rotations,'' Clarendon Press, Oxford, 1964.

\bibitem{Erler}
J.~Erler, {\em Anomaly cancellation in six dimensions}, J. Math. Phys. {\bf
35} (1994) 1819--1833, {\urlfont hep-th/9304104}.

\bibitem{F}W. Fulton, ``Intersection Theory,''  Ergebn. Math.  Grenzegeb.
(3) {\bf 2}, Springer-Verlag, Berlin, 1984.

\bibitem{G}
A.~Grassi, {\em Divisors on elliptic Calabi--Yau 4-folds and the
superpotential in F-theory, I}, J. Geom. Phys. {\bf
28} (1998) 289--319, {\urlfont alg-geom/9704008}.

\bibitem{AD}
A.~Grassi and D.~R.~Morrison {\em Anomalies and the
Euler characteristic of elliptic Calabi--Yau threefolds}, in preparation.

\bibitem{GS}
M.~B.~Green and J.~H.~Schwarz,
{\em Anomaly cancellations in supersymmetric {$D = 10$} gauge theory and
superstring theory},
Phys. Lett. B {\bf 149} (1984) 117--122.

\bibitem{IMS}
K.~Intriligator, D. R. Morrison, and N.~Seiberg, {\em Five-dimensional
supersymmetric gauge
  theories and degenerations of {C}alabi--{Y}au spaces}, Nucl. Phys. B {\bf
  497} (1997) 56--100, {\urlfont hep-th/9702198}.

\bibitem{IN} Y. Ito and H. Nakajima, {\em McKay correspondence and Hilbert
schemes in dimension three}, Topology, to appear, {\urlfont
alg-geom/9803120}.

\bibitem{IR} Y. Ito and M. Reid, {\em The McKay correspondence for finite
subgroups of $SL(3, \mathbb C$)}, in: ``Higher Dimensional Varieties (Trento
1994),'' de~Gruyter, Berlin, 1996, pp. 221--240, {\urlfont
alg-geom/9411010}.

\bibitem{KM} S. Katz and D. R. Morrison, {\em Gorenstein threefold
singularities with small resolutions via invariant theory for Weyl groups},
J. Algebraic Geom. {\bf 1}
     (1992) 449-530, {\urlfont alg-geom/9202002}.

\bibitem{KMP} S. Katz, D. R. Morrison, and M. R. Plesser, {\em
Enhanced gauge symmetry in type {II} string
  theory}, Nucl. Phys. B {\bf 477} (1996) 105--140, {\urlfont
  hep-th/9601108}.

\bibitem{KV} S. Katz and C. Vafa, {\em Matter from geometry},
Nucl. Phys. B {\bf 497} (1997) 146--154, {\urlfont hep-th/9606086}.

\bibitem{kod} K. Kodaira, {\em On compact analytic surfaces, II, III},
Ann. of Math. {\bf 77} (1963) 563--626, {\bf 78} (1963) 1--40.

\bibitem{MP} J. W. G. ~McKay and  J. ~Patera
``Tables of Dimensions, Indices, and Branching Rules for Representations of
Simple Lie Algebras,'' M. Dekker, New York, 1981.

\bibitem{Mi} R. Miranda, {\em Smooth models for elliptic threefolds,} in:
``The Birational Geometry of Degenerations,''
 Progr. Math. {\bf 29}, Birkh\"auser, Boston, 1983, pp. 85--133.

\bibitem{MVII} D. R. Morrison and C. Vafa, {\em
Compactifications of {F}-theory on {C}alabi--{Y}au threefolds, II},
Nucl. Phys. B {\bf 476} (1996) 437--469, {\urlfont
  hep-th/9603161}.

\bibitem{Na1} N. Nakayama, {\em On Weierstrass models}, in: ``Algebraic
Geometry and Commutative Algebra,'' Vol. II, Kinokuniya, Tokyo, 1988,
pp. 405--431.

\bibitem{Na2}
N. Nakayama, {\em Elliptic fibrations over surfaces, I}, in: ``Algebraic
Geometry and Analytic Geometry (Tokyo, 1990),'' Springer, Tokyo, 1991, pp.
126--137.

\bibitem{R} M.~Reid, {\em McKay correspondence}, {\urlfont
  alg-geom/9702016}.

\bibitem{Sadov}
V.~Sadov, {\em Generalized Green--Schwarz mechanism in $F$ theory},
Phys. Lett. B {\bf 388} (1996) 45--50, {\urlfont hep-th/9606008}.

\bibitem{Sagnotti}
A. Sagnotti, {\em A note on the Green--Schwarz mechanism in open-string
theories}, Phys. Lett. B {\bf 294} (1992) 196--203, {\urlfont
hep-th/9210127}.

\bibitem{Schwarz}
J.~H.~Schwarz, {\em Anomaly-free supersymmetric models in six dimensions},
Phys. Lett. B {\bf 371} (1996) 223--230,
{\urlfont hep-th/9512053}.

\bibitem{SVW}
S.~Sethi, C.~Vafa, and E.~Witten, {\em Constraints on low-dimensional string
  compactifications}, Nucl. Phys. B {\bf 480} (1996) 213--224, {\urlfont
  hep-th/9606122}.

\bibitem{witten:MF}
E.~Witten,
\newblock {\em Phase transitions in M-theory and F-theory},
\newblock Nucl. Phys. {\bf B471} (1996) 195--216, {\urlfont hep-th/9603150}.

\end{thebibliography}
\end{document}